\documentclass[11pt,a4paper, reqno]{amsart}
\usepackage{mathrsfs}

\usepackage{amsmath,amsfonts,verbatim}
\usepackage{latexsym}
\usepackage{amssymb,leftidx}
\usepackage{extarrows}
\usepackage{overpic}
\usepackage{color}
\usepackage{epsfig}
\usepackage{subfigure}
\usepackage{tikz}
\usepackage{bbm}
\usepackage{tikz, caption}
\usepackage[font=small,labelfont=bf]{caption}

\usepackage[colorlinks=true,backref=page]{hyperref}
\hypersetup{urlcolor=blue, citecolor=red}

\usepackage{amssymb}
\usepackage{mathrsfs}
\usepackage{amscd}
\usepackage{cite}

\newcommand{\R}{\mathbb{R}}
\newcommand{\C}{\mathbb{C}}
\newcommand{\Sph}{\mathbb{S}}
\newcommand{\Fk}{\mathcal{F}_{\kappa}}

\newcommand{\Lk}{\mathcal{L}_{\kappa}}
\newcommand{\Dk}{\Delta_{\kappa}}
\newcommand{\wk}{\omega_{\kappa}}
\newcommand{\gk}{\gamma_{\kappa}}
\newcommand{\dd}{\,\mathrm{d}}
\newcommand{\Lp}[1]{L^{#1}_{\kappa}(\R^d)}
\newcommand{\norm}[2]{\left\lVert #1\right\rVert_{#2}}
\newcommand{\abs}[1]{\left\lvert #1\right\rvert}

\usepackage{graphicx}
\usepackage{float}

\numberwithin{equation}{section}
\newtheorem{theorem}{Theorem}[section]
\newtheorem{lemma}[theorem]{Lemma}
\newtheorem{proposition}[theorem]{Proposition}

\newtheorem{remark}[theorem]{Remark}

\begin{document}
\title[Uniform Resolvent Estimates for the Dunkl Laplacian ]{Uniform Resolvent Estimates for the Dunkl Laplacian}

\author{Xiaofen Gao}
\address{School of Mathematics and Statistic, Zhengzhou University, Zhengzhou 450001, China}
\email{gaoxiaofen@zzu.edu.cn}

\author{Manli Song}
\address{School of Mathematics and Statistics, Northwestern Polytechnical University, Xi'an, Shaanxi 710129, China}
\email{mlsong@nwpu.edu.cn}

\author{Junyong Zhang}
\address{Department of Mathematics, Beijing Institute of Technology, Beijing 100081}
\email{zhang\_junyong@bit.edu.cn}

\maketitle 

\begin{abstract}Let \(\Lk=-\Dk\) denote the nonnegative Dunkl Laplacian associated with a finite reflection group and a nonnegative multiplicity function \(\kappa\), and let
$
N=d+2\gk\geq 3
$
be the corresponding homogeneous dimension. We establish uniform \(L_\kappa^p(\mathbb{R}^d)\)-\(L_\kappa^q(\mathbb{R}^d)\) resolvent estimates for \(\Lk\) whenever
\[
\frac{2}{N+1}
\leq \frac{1}{p}-\frac{1}{q}
\leq \frac{2}{N},
\quad
\frac{2N}{N+3}<p<\frac{2N}{N+1},
\quad
\frac{2N}{N-1}<q<\frac{2N}{N-3},
\]
where the last upper bound is interpreted as \(q<\infty\) when \(N=3\). Remarkably, this admissible range coincides with the corresponding range for the Euclidean Laplacian \(-\Delta\) in dimension \(N\).

\end{abstract}

\begin{center}
 \begin{minipage}{120mm}
 \end{minipage}
\end{center}

\section{Introduction}



Uniform resolvent estimates for elliptic operators form a classical theme at
the intersection of harmonic analysis, spectral theory, and partial
differential equations
\cite{KRS1987,FrankSabin2017,BoucletMizutani2018}. Beyond their role as
scale-invariant Sobolev inequalities with a spectral parameter, such estimates
underpin limiting absorption principles and have important applications to
unique continuation, spectral enclosure and stability questions for
non-self-adjoint perturbations, scattering theory, and dispersive estimates
\cite{Agmon1975,KRS1987,Frank2011,BoucletMizutani2018,FrankSabin2017}. For the
Euclidean Laplacian, Kenig, Ruiz, and Sogge \cite{KRS1987} established the foundational
uniform Sobolev inequalities on the Sobolev line, whereas
Guti\'errez  \cite{Gutierrez2004} obtained free-resolvent estimates in a larger admissible
$L^p$--$L^q$ region by exploiting both the oscillatory structure of the kernel
and Fourier restriction estimates. Indeed, they have proved for $d\geq3$, there holds
\begin{equation}\label{eq:KRS}
 \bigl\|(-\Delta-z)^{-1}f\bigr\|_{L^q(\R^d)}
 \leq
 C|z|^{\frac{d}{2}(\frac1p-\frac1q)-1}
 \|f\|_{L^p(\R^d)},\quad z\in\C\setminus[0,\infty),
\end{equation}
provided that $(1/p,1/q)$ is contained in the quadrilateral
\begin{equation*}\label{pq2}
	\Big\{(x,y): \frac 2{d+1}\le x-y \leq\frac2d, \, \frac{d+1}{2d}<x<\frac{d+3}{2d}, \,\frac{ d-3}{2d}< y<\frac{d-1}{2d} \Big\}.	
\end{equation*}	
More precisely, for a
boundary spectral parameter $z=\lambda\pm i0$ with $\lambda>0$, the multiplier
$(|\xi|^2-z)^{-1}$ is singular on the characteristic sphere
$|\xi|=\sqrt{\lambda}$. Restriction extension theory and oscillatory integral
analysis therefore enter the resolvent problem in an essential way. This
places uniform resolvent estimates within the same broad circle of ideas as
the Fourier restriction conjecture and its connections with the
Bochner--Riesz, Kakeya, and local smoothing conjectures, see Tao's survey
\cite{Tao2004} for an overview.

\vspace{0.2cm}

The aim of the present paper is to establish uniform
$L^p_\kappa(\R^d)\to L^q_\kappa(\R^d)$ resolvent estimates for the
nonnegative Dunkl Laplacian as above in \eqref{eq:KRS}.  We first recall the underlying
geometric and analytic framework.  Let $R\subset\R^d$ be a reduced root
system, let $G$ be the associated finite reflection group, and let
$\kappa\colon R\to[0,\infty)$ be a $G$-invariant multiplicity function.
After fixing a positive subsystem $R_+$, set
\[
 \wk(x)=\prod_{\alpha\in R_+}
 |\langle x,\alpha\rangle|^{2\kappa(\alpha)},
 \quad
 \gk=\sum_{\alpha\in R_+}\kappa(\alpha),
 \quad
 N=d+2\gk.
\]
Since $\wk(tx)=t^{2\gk}\wk(x)$ for $t>0$, the measure
$d\mu_\kappa(x)=\wk(x)\,dx$ obeys the scaling identity
\begin{equation*}\label{eq:intro-measure-scaling}
 \int_{\R^d}F(tx)\,d\mu_\kappa(x)
 =t^{-N}\int_{\R^d}F(y)\,d\mu_\kappa(y),
 \quad t>0.
\end{equation*}
Accordingly, $N$ is the homogeneous dimension of the Dunkl
structure.  We write
$L^r_\kappa(\R^d)=L^r(\R^d,d\mu_\kappa)$. In particular, if
$D_t f(x)=f(tx)$, then
\begin{equation}\label{eq:intro-lp-scaling}
 \|D_t f\|_{L^r_\kappa}=t^{-N/r}\|f\|_{L^r_\kappa},
 \quad 1\le r\le\infty.
\end{equation}

For $\xi\in\R^d$, the Dunkl operator associated with $\xi$ is the
differential--difference operator
\begin{equation*}\label{eq:intro-dunkl-operator}
 T_\xi f(x)=\partial_\xi f(x)
 +\sum_{\alpha\in R_+}\kappa(\alpha)\langle\alpha,\xi\rangle
 \frac{f(x)-f(\sigma_\alpha x)}{\langle\alpha,x\rangle},
\end{equation*}
where $\sigma_\alpha$ denotes reflection in the hyperplane orthogonal to
$\alpha$.  These mutually commuting operators were introduced by Dunkl
\cite{dunkl_1989}.  Writing $T_j=T_{e_j}$, one defines
\[
 \Dk=\sum_{j=1}^dT_j^2,
 \quad
 \Lk=-\Dk.
\]
The homogeneity of the Dunkl operators gives
\begin{equation}\label{eq:intro-operator-scaling}
 \Lk D_t=t^2D_t\Lk.
\end{equation}
The Dunkl transform, developed systematically by de~Jeu
\cite{dejeu_1993}, diagonalizes $\Lk$
\begin{equation*}\label{eq:intro-dunkl-diagonalization}
 \Fk(\Lk f)(\xi)=|\xi|^2\Fk f(\xi).
\end{equation*}
When $\kappa\equiv0$, this reduces to the usual Fourier-transform
diagonalization of the Euclidean Laplacian.  Consequently, for
$z\in\C\setminus[0,\infty)$,
\begin{equation*}\label{eq:intro-resolvent-multiplier}
 (\Lk-z)^{-1}f
 =\Fk^{-1}\!\left(
   \frac{\Fk f(\xi)}{|\xi|^2-z}
 \right).
\end{equation*}
Moreover, \eqref{eq:intro-lp-scaling} and
\eqref{eq:intro-operator-scaling} show that any scaling invariant uniform
resolvent estimate must have the form
\begin{equation*}\label{eq:intro-resolvent-scaling}
 \| (\Lk-z)^{-1}f\|_{L^q_\kappa}
 \le C |z|^{\frac N2(\frac1p-\frac1q)-1}
 \|f\|_{L^p_\kappa},
\end{equation*}
for the admissible pairs $(p,q)$ and spectral parameters under
consideration.  Thus the homogeneous dimension $N$, rather than the topological
dimension $d$, governs both the admissible exponent geometry and the
spectral scaling of the Dunkl resolvent.
\vspace{0.2cm}

In the study of uniform resolvent estimates, the Euclidean theory relies on several tools that are closely adapted to the classical Schr\"odinger structure, including resolvent identities, multiplier and commutator methods, weighted estimates, and Fourier restriction theory. For example, see \cite{KRS1987,Gutierrez2004} for the free Schr\"odinger operator, \cite{BurqEtAl2003,BoucletMizutani2018,Mizutani2020} for Schr\"odinger operators with electric perturbations, and \cite{CueninKenig2017,FanelliZhangZheng2023} for magnetic perturbations. In particular, Fanelli, Zhang, and Zheng \cite{FanelliZhangZheng2023} obtained $L^p$--$L^q$ uniform resolvent estimates for two-dimensional scaling-critical magnetic Schr\"odinger operators, including the Aharonov--Bohm model, by exploiting the explicit structure of the magnetic resolvent kernel together with dyadic localization, oscillatory integral estimates, and $TT^*$ arguments. Such results illustrate that, near the continuous spectrum, uniform resolvent estimates are governed not only by the local singularity of the kernel but also by the oscillatory behavior associated with the characteristic energy surface.

\vspace{0.2cm}

Nevertheless, a direct transplantation of the Euclidean convolution argument to the Dunkl setting is not available, mainly because the general $L^p$ theory for arbitrary Dunkl translations and convolutions needed for such an approach is not presently available. The decisive observation in the present work is that the resolvent multiplier
$
1/(|\xi|^2-z)
$
is radial. Consequently, its inverse Dunkl transform is radial as well, and the radial Dunkl transform reduces to a Fourier--Bessel transform of order
$
\nu=\frac{N}{2}-1
$
with $N=d+2\gamma_\kappa.$
This reduction transfers the essential analysis from the general Dunkl translation structure to a one-dimensional Bessel integral governed by the homogeneous dimension $N$. In particular, the resolvent kernel exhibits the same two basic regimes as the Euclidean free resolvent in dimension $N$. Near the origin it has the potential-type singularity $r^{2-N}$, whereas in the far field its incoming and outgoing components have the form
$
e^{\pm ir}r^{-(N-1)/2}a_\pm(r),
$
with $a_\pm$ satisfying symbol-type estimates.
\vspace{0.2cm}

These two regimes are controlled by different analytic mechanisms. For the local part, R\"osler's positive radial product formula \cite{Rosler2003} provides the positivity and convolution properties needed to compare the resolvent kernel with a Dunkl--Riesz potential, while the Hardy--Littlewood--Sobolev inequality of Hassani, Mustapha, and Sifi \cite{HassaniMustaphaSifi2009} yields the corresponding fractional-integration estimates on the Sobolev line. For the oscillatory far field part, dyadic localization at the spatial scale $r\sim 2^j$ leads, after the Fourier--Bessel transform, to multipliers concentrated in a $2^{-j}$-neighborhood of the characteristic sphere $|\xi|=1$, with radial $L^2$ size of order $2^{j/2}$. Combining this concentration phenomenon with the Fourier--Dunkl sphere restriction theorem \cite{SenapatiEtAl2024} yields the Stein--Tomas type dyadic estimates required for the limiting absorption argument. Thus the proof is built on two complementary mechanisms dictated by the homogeneous dimension $N$, that is, the local Riesz-potential structure and the far field restriction-driven oscillation. The principal contribution consists in linking these mechanisms via the radial structure of the Dunkl resolvent, in a context in which the standard Euclidean translation argument is not available.


\vspace{0.2cm}
We are now in a position to state our main theorem.

\begin{theorem}[Uniform resolvent estimate]\label{thm:target}
Assume that $N=d+2\gk\geq3$. Let $p$ and $q$ satisfy
\[
 \frac{2}{N+1}
 \leq \frac{1}{p}-\frac{1}{q}
 \leq \frac{2}{N},
 \quad
 \frac{2N}{N+3}<p<\frac{2N}{N+1},
 \quad
 \frac{2N}{N-1}<q<\frac{2N}{N-3}.
\]
Then there exists a constant $C=C(\kappa,N,p,q)>0$, independent of
$z\in\C\setminus[0,\infty)$, such that
\begin{equation}\label{eq:target}
 \bigl\|(\Lk-z)^{-1}f\bigr\|_{L^q_\kappa(\R^d)}
 \leq
 C|z|^{\frac{N}{2}(\frac1p-\frac1q)-1}
 \|f\|_{L^p_\kappa(\R^d)}, \quad \forall f\in L^p_\kappa(\R^d).
\end{equation}

\end{theorem}

\begin{remark}
As shown in Figure~1, the ranges of \(p\) and \(q\) agree with the validity range of the uniform resolvent estimates for \(-\Delta\) in dimension \(N\).
\begin{figure}[htbp]
\centering

\def\dimpar{5}
\def\axislen{7}

\pgfmathsetmacro{\xmin}{(\dimpar+1)/(2*\dimpar)}
\pgfmathsetmacro{\xmax}{(\dimpar+3)/(2*\dimpar)}
\pgfmathsetmacro{\ymin}{(\dimpar-3)/(2*\dimpar)}
\pgfmathsetmacro{\ymax}{(\dimpar-1)/(2*\dimpar)}

\pgfmathsetmacro{\yB}{(\dimpar-1)^2/
  (2*\dimpar*(\dimpar+1))}
\pgfmathsetmacro{\xBp}{(\dimpar^2+4*\dimpar-1)/
  (2*\dimpar*(\dimpar+1))}

\pgfmathsetmacro{\stripLower}{2/(\dimpar+1)}
\pgfmathsetmacro{\stripUpper}{2/\dimpar}

\pgfmathsetmacro{\xE}{(\dimpar^2+\dimpar+2)/(2*\dimpar^2)}
\pgfmathsetmacro{\yE}{((\dimpar-1)*(\dimpar-2))/(2*\dimpar^2)}

\pgfmathsetmacro{\xEp}{(\dimpar^2+3*\dimpar-2)/(2*\dimpar^2)}
\pgfmathsetmacro{\yEp}{(\dimpar^2-\dimpar-2)/(2*\dimpar^2)}

\begin{tikzpicture}[line cap=round,line join=round]

\draw (0,0) rectangle (\axislen,\axislen);
\draw[->] (0,0) -- ({\axislen+0.7},0)
  node[below] {$1/p$};
\draw[->] (0,0) -- (0,{\axislen+0.7})
  node[left] {$1/q$};

\node[below left] at (0,0) {$0$};
\node[below] at (\axislen,0) {$1$};
\node[left] at (0,\axislen) {$1$};

\coordinate (A)  at
  ({\axislen*\xmin},{\axislen*\ymin});
\coordinate (B)  at
  ({\axislen*\xmin},{\axislen*\yB});
\coordinate (Bp) at
  ({\axislen*\xBp},{\axislen*\ymax});
\coordinate (Ap) at
  ({\axislen*\xmax},{\axislen*\ymax});
\coordinate (E)  at
  ({\axislen*\xE},{\axislen*\yE});
\coordinate (Ep) at
  ({\axislen*\xEp},{\axislen*\yEp});

\filldraw[
  fill=gray!40,
  draw=none
]
  (A) -- (Ap) -- (Bp) -- (B) -- cycle;

\draw[dashed,gray]
  ({\axislen*\xmin},0) --
  ({\axislen*\xmin},\axislen);

\draw[dashed,gray]
  ({\axislen*\xmax},0) --
  ({\axislen*\xmax},\axislen);

\draw[dashed,gray]
  (0,{\axislen*\ymin}) --
  (\axislen,{\axislen*\ymin});

\draw[dashed,gray]
  (0,{\axislen*\ymax}) --
  (\axislen,{\axislen*\ymax});



\draw[dashed,thick] (A) -- (B);
\draw[thick]        (B) -- (Bp);
\draw[dashed,thick] (Bp) -- (Ap);
\draw[thick]        (Ap) -- (A);

\foreach \P in {A,B,Bp,Ap}
  \draw[fill=white,thick] (\P) circle (2.2pt);

\node[below right] at (E)  {$E$};
\node[below right] at (Ep) {$E^\prime$};

\node[below left]  at (A)  {$A$};
\node[above left]  at (B)  {$B$};
\node[above]       at (Bp) {$B'$};
\node[above right] at (Ap) {$A'$};

\node[below] at ({\axislen*\xmin},0)
  {$\displaystyle\frac{N+1}{2N}$};

\node[below] at ({\axislen*\xmax},0)
  {$\displaystyle\frac{N+3}{2N}$};

\node[left] at (0,{\axislen*\ymin})
  {$\displaystyle\frac{N-3}{2N}$};

\node[left] at (0,{\axislen*\ymax})
  {$\displaystyle\frac{N-1}{2N}$};

\end{tikzpicture}
\label{figure1}
\caption{Here $A=\left(\frac{N+1}{2N},\frac{N-3}{2N}\right)$, $B=\left(
\frac{N+1}{2N},
\frac{(N-1)^2}{2N(N+1)}
\right)$, $E=\left(
\frac{N^2+N+2}{2N^2},
\frac{(N-1)(N-2)}{2N^2}
\right)$, and $A',B',E'$ are dual points of $A,B,E$, respectively.}
\end{figure}
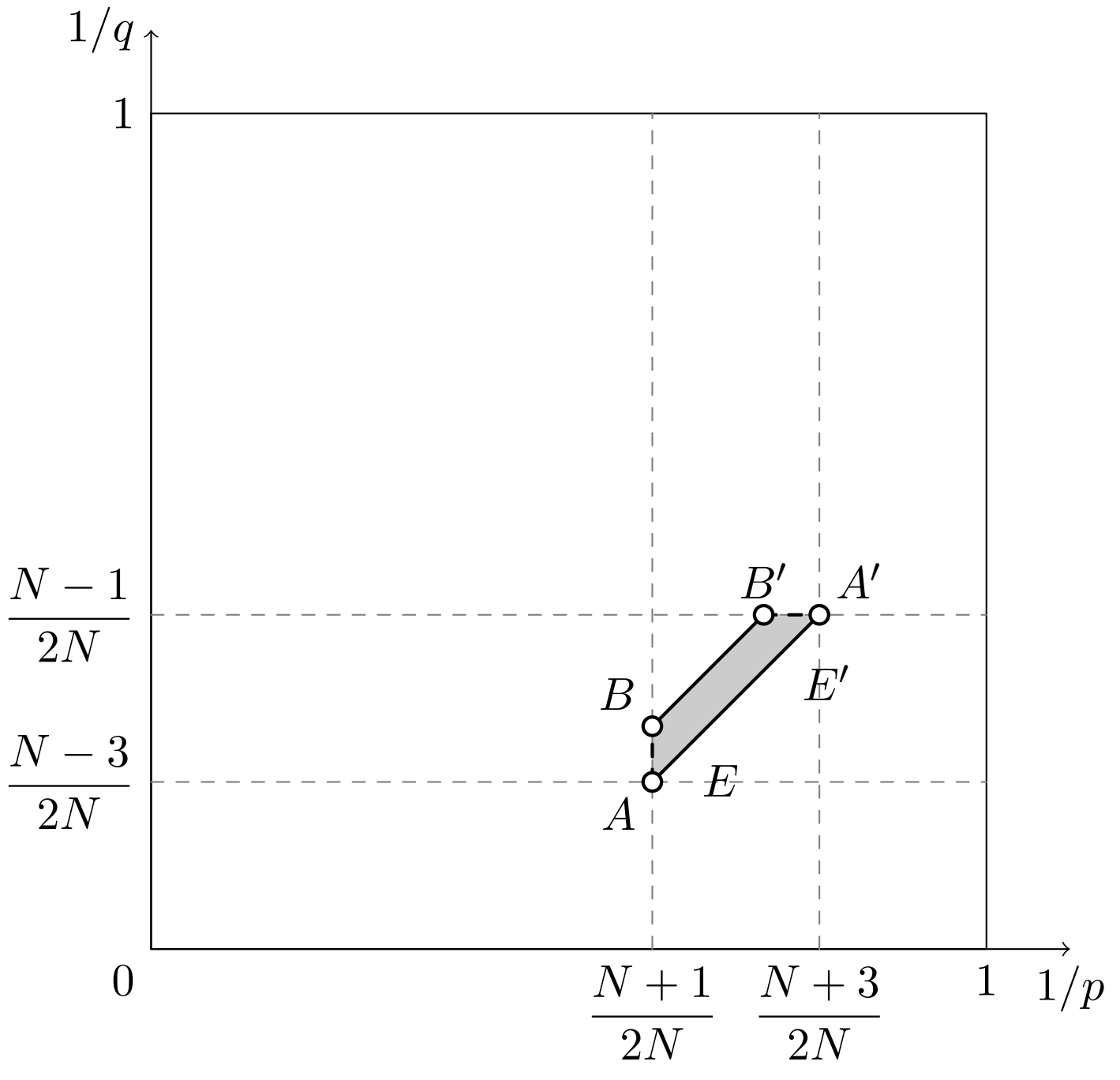
\end{remark}

\begin{remark}\label{rem:scaling}
The power of $|z|$ in \eqref{eq:target} is forced by the homogeneity of $\Lk$
and $d\mu_\kappa$. The upper boundary
$1/p-1/q=2/N$ is the Sobolev line dictated by the local kernel singularity
$r^{2-N}$, whereas the lower boundary
$1/p-1/q=2/(N+1)$ reflects the sphere restriction estimate. Although $\Lk$ is
nonnegative and self-adjoint on $L^2_\kappa(\R^d)$, the off-diagonal
$L^p_\kappa\to L^q_\kappa$ estimate is not a consequence of the $L^2$
spectral theorem, it requires pointwise kernel estimates and oscillatory
analysis.
\end{remark}

We briefly describe the proof.  By scaling,
it suffices to consider $|z|=1$ with $z\neq1$.  Away from the positive spectrum $[0,\infty)$, the
resolvent kernel is expressed by
$$
  K_z(r)
  =c_{\kappa,N}
    \left(\frac{w}{r}\right)^{\nu}K_{\nu}(wr),
  \quad w=\sqrt{-z},\; \operatorname{Re}w>0,
$$
and satisfies
$$
 |K_z(r)|
 \leq C_{\varepsilon,\kappa,N}
 \left(
 r^{2-N}\mathbf{1}_{(0,1]}(r)
 +r^{-(N-1)/2}e^{-c_\varepsilon r}\mathbf{1}_{[1,\infty)}(r)
 \right).
$$
Here $K_\nu$ denotes the modified Bessel function of second kind.
Radial Young inequalities yield the subcritical estimates, while on the
Sobolev line the local singularity is controlled by the Dunkl--Riesz
potential of order two and the Hardy--Littlewood--Sobolev theorem
\cite{HassaniMustaphaSifi2009}. Near the positive spectrum, the incoming and outgoing Hankel kernels are
split into local and oscillatory parts.  The former is handled by the same
potential estimates since it behaves like $r^{2-N}$ .  For the latter, a dyadic decomposition at scale
$R=2^j$ gives
\[
 k_j(r)=\beta(r/R)e^{\pm ir}r^{-(N-1)/2}a_\pm(r),
\]
where $a_\pm$ is a symbol of order zero.  Writing $m_j=\Fk(k_j)$, the
Fourier--Bessel reduction yields, for every $M\geq1$,
\begin{equation}\label{eq:multiplier-concentration}
 |m_j(\rho)|
 \leq C_MR\bigl(1+R|\rho-1|\bigr)^{-M},
 \quad \rho\geq0.
\end{equation}
Thus $m_j$ is concentrated in an $R^{-1}$-neighborhood of the characteristic
sphere and has radial $L^2$ norm $O(R^{1/2})$. Combining \eqref{eq:multiplier-concentration} with the Fourier--Dunkl sphere
restriction theorem \cite{SenapatiEtAl2024} gives
\[
 \|f\ast_D k_j\|_{L^{q_0}_\kappa}
 \lesssim R^{1/2}\|f\|_{L^2_\kappa},
 \quad q_0=\frac{2(N+1)}{N-1}.
\]
Together with the $L^1_\kappa\to L^\infty_\kappa$ kernel bound, interpolation,
duality, and dyadic summation produce the asserted exponent region.  At the
two limiting vertices, restricted weak-type estimates are first established
and then interpolated to obtain the strong-type interior bounds.

\vspace{0.2cm}

The paper is organized as follows. Section~2 contains the required facts about
Dunkl operators, the Dunkl transform, radial convolution, Bessel functions,
and Dunkl--Riesz potentials. Section~3 gives the proof of  our main Theorem \ref{thm:target}. The proof of Theorem \ref {theo:limits} and Lemma \ref{lem:dyadic} are  given in Sections 4 and 5, respectively.

\subsection*{Notations}

\begin{itemize}
    \item \(A \lesssim B\) means that \(A \leqslant CB\), and \(A \sim B\) stands for \(C_1B \leq A \leq C_2B\), where \(C\), \(C_1\), \(C_2\) denote positive constants  that may depend on parameters such as \(p\), \(q\), \(\kappa\), and \(N\).
    \item \(\Delta_\kappa\) denotes Dunkl--Laplacian on \(\mathbb{R}^n\).
    \item $L_\kappa^p(\R^d)$ denotes Dunkl Sobolev space.
    \item $\nu=\nu(N)=\frac{N}{2}-1.$
    \item For any $\mu>-\frac{1}{2}$, \(J_\mu\) denotes the Bessel function of order \(\mu\),
\begin{equation}\label{Bessel}
J_\mu(r)=\frac{(r/2)^\mu}{\Gamma(\mu+1/2)\pi^{1/2}}\int_{-1}^1e^{ir\tau}(1-\tau^2)^{\mu-1/2}d\tau,
\end{equation}
and $\mathfrak j_\mu$ denote the normalized Bessel function
\begin{equation}\label{eq:jnorm}
\mathfrak j_\mu(r)=2^\mu\Gamma(\mu+1)r^{-\mu}J_\mu(r),
\quad
\mathfrak j_\mu(0)=1.
\end{equation}
\item If $f$ is a radial function on $\mathbb{R}^d$, i.e., $f(x)=\tilde{f}(|x|)$, we abuse notation and write $f(r)$ for $\tilde{f}(r)$ for convenience.
\item For any $\mu\in\mathbb{R}$,  $\mathcal{K}_{\mu}$ denotes the modified Bessel function of the second kind
\begin{equation}\label{second}
\mathcal{K}_{\mu}(z)
  =\int_{0}^{\infty}e^{-z\cosh \tau}\cosh(\mu \tau)\,\dd \tau,
  \quad \operatorname{Re}z>0.
  \end{equation}
\end{itemize}


\section{Preliminaries}

In this section, we recall some basic definitions and important properties of the Dunkl
operators and briefly overview the related harmonic analysis to make the paper self-contained.
A complete account of harmonic analysis related to Dunkl operators can be found in
\cite{LMS2026,bensaid_kobayashi_orsted_2012,bensaid_nandakumaran_ratnakumar_2011,
mejjaoli_2009_dispersion,mondal_song_2025_orthonormal,rosler_1999_positivity,thangavelu_xu_2007_riesz}.

\subsection{Dunkl Operators}

The basic ingredient in the theory of Dunkl operators are root systems and finite reflection
groups. We start this section by the definition of root system.

Let $\langle\cdot,\cdot\rangle$ denotes the standard Euclidean scalar product in $\R^d$.
For $x\in\R^d$, we denote $|x|$ as $|x|=\langle x,x\rangle^{1/2}$. For
$\alpha\in\R^d\setminus\{0\}$, we denote $r_\alpha$ as the reflection with respect to
the hyperplane $\langle\alpha\rangle^\perp$, orthogonal to $\alpha$ and is defined by
\[
  \sigma_\alpha(x):=x-2\frac{\langle\alpha,x\rangle}{|\alpha|^2}\alpha,
  \qquad x\in\R^d.
\]
A finite set $R$ in $\R^d\setminus\{0\}$ is said to be a root system if the following holds:
\begin{enumerate}
  \item $\sigma_\alpha(R)=R$ for all $\alpha\in R$,
  \item $R\cap\mathbb{R}\alpha=\{\pm\alpha\}$ for all $\alpha\in R$.
\end{enumerate}

For a given root system $R$, the subgroup $G\subset O(n,\mathbb{R})$ generated by the
reflections $\{\sigma_\alpha\mid\alpha\in R\}$ is called the finite Coxeter group associated with
$R$. The dimension of $\operatorname{span}R$ is called the rank of $R$. For a detailed study
on the theory of finite reflection groups, we refer to \cite{humphreys_1990}. Let
\[
  R_+:=\{\alpha\in R:\langle\alpha,\beta\rangle>0\}
\]
for some $\beta\in\R^d\setminus\bigcup_{\alpha\in R}\langle\alpha\rangle^\perp$, be a
fix positive root system.

Some typical examples of such a system are the Weyl groups such as the symmetric group
$S_d$ for the type $A_{d-1}$ root system and the hyperoctahedral group for the type $B_d$
root system. In addition, $H_3$, $H_4$ (icosahedral groups) and $I_2(d)$ (symmetry group of
the regular $d$-gon) are also the Coxeter groups.

A multiplicity function for $G$ is a function $\kappa:R\to\mathbb{C}$ which is constant on
$G$-orbits. Setting $\kappa_\alpha:=\kappa(\alpha)$ for $\alpha\in R$, from the definition of
$G$-invariant, we have $\kappa_{g\alpha}=\kappa_\alpha$ for all $g\in G$. We say $\kappa$ is
non-negative if $\kappa_\alpha\geq 0$ for all $\alpha\in R$. The $\mathbb{C}$-vector space of
non-negative multiplicity functions on $R$ is denoted by $\mathcal{K}^+$. Let us denote
$\gamma$ as
\[
  \gamma=\gamma(\kappa):=\sum_{\alpha\in R_+}\kappa(\alpha).
\]
For convenience, we denote $N=d+2\gamma$ throughout the paper.

For $\xi\in\mathbb{C}^d$ and $\kappa\in\mathcal{K}^+$, Dunkl in 1989 introduced a family of
first-order differential--difference operators $T_\xi:=T_\xi(\kappa)$ by
\begin{equation}
  T_\xi(\kappa)f(x)
  :=\partial_\xi f(x)
  +\sum_{\alpha\in R_+}\kappa_\alpha\langle\alpha,\xi\rangle
  \frac{f(x)-f(\sigma_\alpha x)}{\langle\alpha,x\rangle},
  \qquad f\in C^1(\R^d),
  \label{def:dunkl-operator}
\end{equation}
where $\partial_\xi$ denotes the directional derivative corresponding to $\xi$. The operator
$T_\xi$ defined in (2.1), formally known as Dunkl operator and is one of the most important
developments in the theory of special functions associated with root systems \cite{dunkl_1989}. They
commute pairwise and are skew-symmetric with respect to the $G$-invariant measure
$\omega_\kappa(x)\,dx$, where the weight function
\[
  \omega_\kappa(x):=\prod_{\alpha\in R_+}|\langle\alpha,x\rangle|^{2\kappa_\alpha}
\]
is of homogeneous degree $2\gamma$. Thanks to the $G$-invariance of the multiplicity
function, the definition of the Dunkl operator is independent of the choice of the positive
subsystem $R_+$. In \cite{dunkl_1991}, it is shown that for any $\kappa\in\mathcal{K}^+$, there is a unique
linear isomorphism $V_\kappa$ (Dunkl's intertwining operator) on the space
$\mathcal{P}(\R^d)$ of polynomials on $\R^d$ such that
\begin{enumerate}
  \item $V_\kappa\bigl(\mathcal{P}_m(\R^d)\bigr)=\mathcal{P}_m(\R^d)$ for all
  $m\in\mathbb{N}$,
  \item $V_\kappa|_{\mathcal{P}_0(\R^d)}=\mathrm{id}$,
  \item $T_\xi(\kappa)V_\kappa=V_\kappa\partial_\xi$,
\end{enumerate}
where $\mathcal{P}_m(\R^d)$ denotes the space of homogeneous polynomials of degree
$m$.

For any finite reflection group $G$ and any $\kappa\in\mathcal{K}^+$, R\"osler in \cite{rosler_1999_positivity} proved
that there exists a unique positive Radon probability measure $\rho_x^\kappa$ on
$\R^d$ such that
\begin{equation}
  V_\kappa f(x)=\int_{\R^d}f(\xi)\,d\rho_x^\kappa(\xi).
  \label{def:dunkl-kernel}
\end{equation}
The measure $\rho_x^\kappa$ depends on $x\in\R^d$ and its support is contained in
the ball
\[
  B(\|x\|):=\{\xi\in\R^d:\|\xi\|\leq\|x\|\}.
\]
In view of the Laplace type representation (2.2), Dunkl's intertwining operator $V_\kappa$
can be extended to a larger class of function spaces.

Let $\{\xi_1,\xi_2,\ldots,\xi_d\}$ be an orthonormal basis of
$(\R^d,\langle\cdot,\cdot\rangle)$. Then the Dunkl Laplacian operator $\Delta_\kappa$
is defined as
\[
  \Delta_\kappa=\sum_{j=1}^d T_{\xi_j}^2(\kappa).
\]
The definition of $\Delta_\kappa$ is independent of the choice of the orthonormal basis of
$\R^d$. In fact, one can see that the operator $\Delta_\kappa$ also can be expressed as
\[
  \Delta_\kappa f(x)
  =\Delta f(x)
  +\sum_{\alpha\in R_+}\kappa_\alpha
  \left\{
    \frac{2\langle\nabla f(x),\alpha\rangle}{\langle\alpha,x\rangle}
    -|\alpha|^2\frac{f(x)-f(\sigma_\alpha x)}{\langle\alpha,x\rangle^2}
  \right\},
  \qquad f\in C^1(\R^d),
\]
where $\nabla$ and $\Delta$ are the usual gradiant and Laplacian operator on $\R^d$,
respectively. Observe that, for $\kappa\equiv 0$, the Dunkl Laplacian operator $\Delta_\kappa$
reduces to the classical Euclidean Laplacian $\Delta$.

For every $y\in\R^d$, the system
\[
  \begin{cases}
    T_\xi u(x,y)=\langle y,\xi\rangle u(x,y), & \xi\in\R^d,\\
    u(0,y)=1,
  \end{cases}
\]
admits a unique analytic solution on $\R^d$, which we denote as $E_\kappa(x,y)$ and is
known as the Dunkl kernel in general. The kernel has a unique holomorphic extension to
$\mathbb{C}^d\times\mathbb{C}^d$ and is given by
\[
  E_\kappa(x,y)=V_\kappa\bigl(e^{\langle\cdot,y\rangle}\bigr)(x)
  =\int_{\R^d}e^{\langle\xi,y\rangle}\,d\rho_x^\kappa(\xi),
  \qquad \forall x,y\in\R^d,
\]
where the positive Radon probability measure $\rho_x^\kappa$ depends on
$x\in\R^d$ and its support is contained in the ball
$B(\|x\|):=\{\xi\in\R^d:\|\xi\|\leq\|x\|\}$. When $\kappa\equiv 0$, the Dunkl
kernel $E_\kappa(x,y)$ reduces to the exponential $e^{\langle x,y\rangle}$. For
$x,y\in\R^d$ and $z,w\in\mathbb{C}^n$, the Dunkl kernel satisfies the following
properties:
\begin{enumerate}
  \item[(i)] $E_\kappa(z,w)=E_\kappa(w,z)$ and $E_\kappa(\lambda z,w)=E_\kappa(z,\lambda w)$, $\lambda\in\mathbb{C}$.
  \item[(ii)] For any $h \in\mathbb{N}^d$, we have
  \[
    |D_z^h E_\kappa(x,z)|\leq |x|^{|h|}e^{|\operatorname{Re}(z)||x|},
  \]
  where
  \[
    D_z^h=\frac{\partial^{|h|}}{\partial z_1^{h_1}\cdots\partial z_d^{h_d}}
    \quad\text{and}\quad
    |h|=h_1+\cdots+h_d.
  \]
  In particular, for any $x,y\in\R^d$, we have
  \[
    |E_\kappa(x,iy)|\leq 1.
  \]
  \item[(iii)] $E_\kappa(ix,y)=\overline{E_\kappa(-ix,y)}$ and $E_\kappa(gx,gy)=E_\kappa(x,y)$, $g\in G$.
\end{enumerate}

\subsection{The Dunkl Transform}

For $1\leq p<\infty$, let $L_\kappa^p(\R^d)$ be the space of $L^p$-functions on
$\R^d$ with respect to the weight $\omega_\kappa(x)$ with the
$L_\kappa^p(\R^d)$-norm
\[
  \|f\|_{L_\kappa^p(\R^d)}
  =\left(\int_{\R^d}|f(x)|^p \omega_\kappa(x)\,dx\right)^{1/p},
\]
and
\[
  \|f\|_{L_\kappa^\infty(\R^d)}
  =\operatorname*{ess\,sup}_{x\in\R^d}|f(x)|.
\]
The Dunkl transform of $f\in L_\kappa^1(\R^d)$ is defined by
\begin{equation}
  \mathcal{F}_\kappa f(y)
  =\frac{1}{c_\kappa}\int_{\R^d}f(x)E_\kappa(x,-iy)\omega_\kappa(x)\,dx,
  \label{def:dunkl-transform}
\end{equation}
where $c_\kappa$ is the Mehta-type constant defined by
\begin{equation*}
   c_\kappa=\int_{\R^d}e^{-|x|^2/2}\omega_\kappa(x)\,dx.
\end{equation*}
The weighted polar decomposition takes the form
\begin{equation}\label{eq:polar}
\int_{\R^d}F(\xi)\wk(\xi)\dd\xi
=c_\kappa
\int_0^\infty\int_{\Sph^{d-1}}
F(\rho\theta)\omega_\kappa(\theta)\,\dd\sigma(\theta)\,
\rho^{N-1}\dd\rho,
\end{equation}
where $\dd\sigma(\theta)$ is the classical spherical measure.
\vspace{0.2cm}

When $\kappa\equiv 0$, the Dunkl transform coincides with the classical Fourier transform.
The Dunkl transform plays the same role as the Fourier transform in classical Fourier analysis
and enjoys properties similar to those of the classical Fourier transform. With
$N=d+2\gamma$, here we list some basic properties of the Dunkl transform (see \cite{dai_ye_2019,dunkl_1992,dejeu_1993}).
\begin{enumerate}
  \item[(i)] For all $f\in L_\kappa^1(\R^d)$, we have
  \[
    \|\mathcal{F}_\kappa f\|_{L_\kappa^\infty(\R^d)}
    \leq\frac{1}{c_\kappa}\|f\|_{L_\kappa^1(\R^d)}.
  \]
  \item[(ii)] For any function $f$ in the Schwartz space $\mathcal{S}(\R^d)$, we have
  \[
    \mathcal{F}_\kappa(T_\xi f)(y)
    =i\langle\xi,y\rangle\mathcal{F}_\kappa f(y),
    \qquad y,\xi\in\R^d.
  \]
  In particular, it follows that
  \[
    \mathcal{F}_\kappa(\Delta_\kappa f)(y)
    =-|y|^2\mathcal{F}_\kappa f(y),
    \qquad y\in\R^d.
  \]
  \item[(iii)] For all $f\in L_\kappa^1(\R^d)$, we have
  \[
    \mathcal{F}_\kappa\bigl(f(\cdot/\lambda)\bigr)(y)
    =\lambda^N\mathcal{F}_\kappa f(\lambda y),
    \qquad y\in\R^d,
  \]
  for all $\lambda>0$.
  \item[(iv)] The Dunkl transform $\mathcal{F}_\kappa$ is a homeomorphism of the Schwartz
  space $\mathcal{S}(\R^d)$ and its inverse is given by
  $\mathcal{F}_\kappa^{-1}g(x)=\mathcal{F}_\kappa g(-x)$, for all
  $g\in\mathcal{S}(\R^d)$. In addition, for all
  $f\in\mathcal{S}(\R^d)$, it satisfies
  \[
    \int_{\R^d}|f(x)|^2\omega_\kappa(x)\,dx
    =\int_{\R^d}|\mathcal{F}_\kappa f(y)|^2\omega_\kappa(y)\,dy.
  \]
  In particular, the Dunkl transform extends to an isometric isomorphism on
  $L_\kappa^2(\R^d)$. The Dunkl transform can be extended to the space of tempered
  distributions $\mathcal{S}'(\R^d)$ and is also a homeomorphism of
  $\mathcal{S}'(\R^d)$.
  \item[(v)] For every $f\in L_\kappa^1(\R^d)$ such that
  $\mathcal{F}_\kappa f\in L_\kappa^1(\R^d)$, we have the inversion formula
  \[
    f(x)=\frac{1}{c_\kappa}\int_{\R^d}
    \mathcal{F}_\kappa f(y)E_\kappa(ix,y)\omega_\kappa(y)\,dy,
    \qquad \text{a.e. }x\in\R^d.
  \]
  \item[(vi)]  Let $
\nu:=\nu(N)=\frac{N}{2}-1$. If $f$ is a radial function in $L_\kappa^1(\R^d)$ such that
 $f(x)=f(\abs{x})$ , then
\begin{equation}\label{eq:radialDunkl}
\begin{aligned}
    \Fk f(\rho)=\Fk f(\xi)
=&\frac{\rho^{-\nu}}{\Gamma(\nu+1)}\int_0^\infty
f(r) J_\nu(r\rho)r^{\nu+1}\dd r\\
=&\frac{1}{2^\nu\Gamma^2(\nu+1)}\int_0^\infty
f(r)\mathfrak j_\nu(r\rho)r^{N-1}\dd r,
\end{aligned}
\end{equation}
where $\rho=|\xi|$, $\forall \xi\in\mathbb{R}^n$.
Thus the radial Dunkl transform coincides with the Fourier--Bessel transform in homogeneous dimension $N$.

\end{enumerate}

\subsection{Dunkl Convolution Operator}

For given $x\in\R^d$, the Dunkl translation operator $f\mapsto\tau_x f$ is defined
on $\mathcal{S}(\R^d)$ by
\[
  \mathcal{F}_\kappa(\tau_x f)(y)
  =E_\kappa(x,iy)\mathcal{F}_\kappa f(y),
  \qquad y\in\R^d.
\]
At the moment, the explicit formula of the Dunkl translation operator is known only in two
cases. One is when $f$ is a continuous radial function in
$L_\kappa^2(\R^d)$, the Dunkl translation operator is represented by (see \cite{dai_wang_2010,rosler_1999_positivity})
\begin{equation}\label{radial}
  \begin{aligned}
    \tau_x f(y)
    &=V_\kappa\left[
      f\left(
        \sqrt{|y|^2+|x|^2-2\langle y,\cdot\rangle}
      \right)
    \right](x)\\
    &=\int_{\R^d}
      f\left(
        \sqrt{|y|^2+|x|^2-2\langle y,\xi\rangle}
      \right)
      \,d\rho_x^\kappa(\xi).
  \end{aligned}
\end{equation}
The Dunkl translation operator can also be extended to the space of tempered distributions
$\mathcal{S}'(\R^d)$. It turns out to be rather difficult to extend important results
in classical Fourier analysis to the setting of the Dunkl transform in general. One of the
difficulties comes from the fact that the Dunkl translation operator is not positive in general.
In fact, even $L_\kappa^p(\R^d)$-boundedness of $\tau_x$ is not established in general.
On the other hand, when $f$ is a radial function in $L_\kappa^p(\R^d)$,
$1\leq p\leq\infty$, the following holds (see \cite{gorbachev_ivanov_tikhonov_2019,thangavelu_xu_2005_convolution})
\[
  \|\tau_x f\|_{L_\kappa^p(\R^d)}
  \leq\|f\|_{L_\kappa^p(\R^d)}.
\]
Using the Dunkl translation operator, we define the Dunkl convolution of functions
$f,g\in\mathcal{S}(\R^d)$ by
\[
  f*_D g(x)
  =\int_{\R^d}\tau_x f(-y)g(y)\omega_\kappa(y)\,dy,
  \qquad x\in\R^d.
\]

The Dunkl convolution satisfies the following properties (see \cite{thangavelu_xu_2005_convolution}).
\begin{enumerate}
  \item $\mathcal{F}_\kappa(f*_D g)
  =\mathcal{F}_\kappa(f)\mathcal{F}_\kappa(g)$ and $f*_D g=g*_D f$.
  \item Young's inequality: Let $1\leq p,q,r\leq\infty$ such that
  $1+\frac{1}{q}=\frac{1}{p}+\frac{1}{r}$. If $f\in L_\kappa^p(\R^d)$ and $g$ is
  a radial function of $L_\kappa^r(\R^d)$, then
  $f*_D g\in L_\kappa^q(\R^d)$ and we have
  \begin{equation}
    \|f*_D g\|_{L_\kappa^q(\R^d)}
    \leq
    \|f\|_{L_\kappa^p(\R^d)}
    \|g\|_{L_\kappa^r(\R^d)}.
    \label{ine:dunkl-young}
  \end{equation}
\end{enumerate}

\subsection{ Technical Lemmas}

In this subsection, we collect several auxiliary results used in the proof of
the main theorem.  After applying the Dunkl transform and the relevant
Young inequalities, the required estimates reduce to oscillatory
integrals involving Bessel functions.  We therefore recall the Dunkl
Hardy--Littlewood--Sobolev inequality, together with the stationary-phase
estimates and Bessel-function properties needed below.

For $0<s<N$ and $f\in\mathcal S(\R^d)$, define the Dunkl--Riesz potential
of order $s$ by
\begin{equation*}
  I_s^\kappa f(x)
  :=\frac 1{d_\kappa^s}
    \int_{\R^d}
      \tau_x\bigl(\lvert\,\cdot\,\rvert^{s-N}\bigr)(-y)
      f(y)\omega_\kappa(y)\,dy,
\end{equation*}
where $$d_{\kappa}^s=2^{-\frac N2+s}\Gamma(\frac s2)/\Gamma(\frac{N-\alpha}2).$$  The following result is the Dunkl analogue of the classical
Hardy--Littlewood--Sobolev inequality, see
\cite[Theorem~1.1(i)]{HassaniMustaphaSifi2009}.

\begin{lemma}[Dunkl Hardy--Littlewood--Sobolev inequality]
\label{lem:dunkl-HLS}
Let $0<s<N$ and $1<p<N/s$, and define $q$ by
\[
  \frac1q=\frac1p-\frac{s}{N}.
\]
Then $I_s^\kappa$ extends uniquely to a bounded operator from
$L_\kappa^p(\R^d)$ to $L_\kappa^q(\R^d)$, and
\begin{equation*}
  \norm{I_s^\kappa f}{L_\kappa^q(\R^d)}
  \le C_{p,s}
       \norm{f}{L_\kappa^p(\R^d)}.
\end{equation*}
\end{lemma}

Let $j_\mu$ be the normalized Bessel function in \eqref{eq:jnorm}, then the following lemma concerning estimates of $j_\mu$ will also be needed.

\begin{lemma}
\label{lem:normalized-bessel-properties}
Let $\mu>-\frac{1}{2}$. For any integer $\ell\ge0$, we have
\begin{equation}
\label{eq:besselsmall}
  \abs{\partial_r^\ell\mathfrak j_\mu(r)}\le1,
  \quad r\geq0.
\end{equation}
Moreover, there exist
$b_+,b_-\in C^\infty(0,\infty)$ such that, for $r\ge1$,
\begin{equation}
\label{eq:bessellarge}
  \mathfrak j_\mu(r)
  =r^{-(\mu+1/2)}
   \bigl(e^{ir}b_+(r)+e^{-ir}b_-(r)\bigr),
\end{equation}
and, for every integer $\ell\ge0$,
\begin{equation}
\label{eq:bsymbol}
  \abs{\partial_r^\ell b_\pm(r)}
  \le C_{\mu,\ell}r^{-\ell},
  \quad r\ge1.
\end{equation}
\end{lemma}

\begin{proof}
By \eqref{Bessel} and \eqref{eq:jnorm}, we have
\begin{equation}\label{eq:poisson-representation}
  \mathfrak j_\mu(r)
  =c_\mu\int_{-1}^{1}e^{ir\tau}
       (1-\tau^2)^{\mu-1/2}\,d\tau,
  \quad
  c_\mu=\frac{\Gamma(\mu+1)}
              {\pi^{1/2}\Gamma(\mu+1/2)}.
\end{equation}
Since $\mu>-1/2$, differentiation under the integral is justified at
every order, and hence, for any $r\geq0$,
\[
  \partial_r^\ell\mathfrak j_\mu(r)
  =c_\mu\int_{-1}^{1}(i\tau)^\ell e^{ir\tau}
       (1-\tau^2)^{\tau-1/2}\,d\tau.
\]
Because $|\tau|\leq1$, it follows that
\begin{equation*}
  \bigl|\partial_r^\ell\mathfrak j_\mu(r)\bigr|
  \leq c_\mu\int_{-1}^{1}(1-\tau^2)^{\mu-1/2}\,d\tau=\mathfrak j_\mu(0)=1,
  \quad \forall r\geq0,
\end{equation*}
which proves \eqref{eq:besselsmall}.

We next establish the large-argument decomposition. Recall that
\begin{equation*}
J_\mu(r)=\frac{(\frac{r}{2})^\mu}{\pi^{1/2}\Gamma(\mu+\frac{1}{2})} \left[ie^{-ir}\int_0^\infty e^{-rt}(t^2+2it)^{\mu-\frac{1}{2}}dt-ie^{ir}\int_0^\infty e^{-rt}(t^2-2it)^{\mu-\frac{1}{2}}dt\right],
\end{equation*}
see \cite[Appendix B]{Grafakos2008}. Then, we have
\begin{equation*}
  \mathfrak j_\mu(r)
  =r^{-(\mu+1/2)}
   \bigl(e^{ir}b_+(r)+e^{-ir}b_-(r)\bigr),
\end{equation*}
where
\begin{equation*}
    b_\pm(r)=\mp i\frac{\Gamma(\mu+1)r^{\mu+1/2}}{\pi^{1/2}\Gamma(\mu+\frac{1}{2})}\int_0^\infty e^{-rt}(t^2\mp 2it)^{\mu-\frac{1}{2}}dt
\end{equation*}
Notice that
\begin{equation*}
   \left|\frac{d^\ell}{dr^\ell}\Big(\int_0^\infty e^{-rt}(t^2\mp 2it)^{\mu-\frac{1}{2}}dt\Big)\right|\lesssim_\ell r^{-(\mu+\frac{1}{2})-\ell},\quad \ell\geq 0, \;r\geq1,
\end{equation*}
which indicates \eqref{eq:bsymbol} and completes the proof.
\end{proof}

\section{The proof of  Theorem \ref{thm:target}}
In this section, we prove Theorem \ref{thm:target}.
By the scaling invariant of the operator $\Lk$, it suffices to prove \eqref{eq:target} when $|z|=1$ with $z\neq 1$, that is,
\begin{equation}\label{est:res-s}
  \bigl\|(\Lk-z)^{-1}f\bigr\|_{\Lp q}
  \leq C\|f\|_{\Lp p},
  \quad |z|=1, z\neq 1,
\end{equation}
where $p,q$ are in Theorem \ref{thm:target}. Indeed, let $t>0$, define
$$D_tf(x):=f(tx).$$ Then, we have$$\Lk D_t=t^2D_t\Lk.$$
Besides,
$$(\Lk-z)D_t=t^{-2}D_t(\Lk-z/t^2).$$
Moreover,
$$(\Lk-z)^{-1}=t^{-2}D_t(\Lk-z/t^2)^{-1}D_{1/t}.$$
Notice that
$$\norm{D_tf}{L^s_\kappa}=t^{-N/s}\norm{f}{L^s_\kappa},\quad\forall 1\leq s<\infty.$$
Then taking $t=|z|^{1/2},\  \zeta=z/|z|,$ there holds
\begin{equation*}
\begin{split}
\|(\Lk-z)^{-1}f\|_{\Lp{q}}&\leq Ct^{-2-N/q}\norm{D_{1/t}f}{L^p_\kappa}\\
&=C t^{-2-N/q+N/p}\norm{f}{L^p_\kappa}\\
&=C|z|^{\frac N2(\frac1p-\frac1q)-1}\|f\|_{\Lp{p}}.
\end{split}
\end{equation*}
Thus, we see that \eqref{eq:target} is reduced to \eqref{est:res-s}.

To prove \eqref{est:res-s}, let
\[
  \Omega_{\varepsilon}
  :=\left\{z\in\C\setminus[0,\infty):
      |z|=1,\ \operatorname{dist}(z,[0,\infty))\geq\varepsilon
    \right\},
\]
and
 \[
  \Omega_{\varepsilon}^{-}
 : =\left\{z\in\C\setminus[0,\infty):
      |z|=1,\ \operatorname{dist}(z,[0,\infty))<\varepsilon
    \right\},
\]
we consider two case: (i) $z\in \Omega_{\varepsilon}$, and (ii) $z\in \Omega_{\varepsilon}^{-}$, where $0<\varepsilon\ll1$.

\subsection{Case (i): $z\in \Omega_{\varepsilon}$.} In this case, spectral parameters remain a
fixed positive distance from $[0,\infty)$. We only need to show the following Proposition.


\begin{proposition}
\label{prop:off-spectrum-resolvent}Let $z\in\Omega_{\varepsilon}$, and suppose that $(p,q)$ satisfies either
$$1\leq p\leq q\leq\infty,\quad
  0\leq \frac1p-\frac1q<\frac2N,
$$
or $$1<p<\frac N2,\quad 0\leq \frac1p-\frac1q=\frac2N.$$
Then, for every $f\in{\Lp p}$,
\[
  \bigl\|(\Lk-z)^{-1}f\bigr\|_{\Lp q}
  \leq C_{\varepsilon}\|f\|_{\Lp p},
\]
uniformly for $z\in\Omega_{\varepsilon}$.
\end{proposition}
The following kernel estimate serves as the key ingredient in the proof of the Proposition \ref{prop:off-spectrum-resolvent}.
\begin{lemma}[Kernel uniform estimates]
\label{lem:positive-real-resolvent-kernel}
Let $N\geq3$ and $0<\varepsilon\ll1$. For every
$z\in\Omega_{\varepsilon}$, let $K_z$ denote the inverse Dunkl transform of
$(|\xi|^{2}-z)^{-1}$. Then  $K_z$ is radial and, for any $r=|x|>0$, satisfies
\begin{equation}
  K_z(r)
  =C_{\kappa,N}
    \left(\frac{w}{r}\right)^{\nu}\mathcal{K}_{\nu}(wr),
  \quad w=\sqrt{-z},\; \operatorname{Re}w>0,
  \label{eq:positive-real-modified-bessel-kernel}
\end{equation}
where $K_{\nu}$ denotes the modified Bessel function of the second kind.
Moreover,
\begin{equation}
  |K_z(r)|
  \leq C_{\varepsilon,\kappa,N}
  \left(
    r^{2-N}\mathbf 1_{(0,1]}(r)
    +r^{-(N-1)/2}e^{-\varepsilon r/2}
       \mathbf 1_{[1,\infty)}(r)
  \right),
  \label{eq:positive-real-two-region-kernel-bound}
\end{equation}
and consequently
\begin{equation}
  |K_z(r)|
  \leq C_{\varepsilon,\kappa,N}r^{2-N},
  \quad r>0,
  \label{eq:positive-real-global-kernel-bound}
\end{equation}
uniformly for $z\in\Omega_{\varepsilon}$.
\end{lemma}
We defer the proof of Lemma  \ref{lem:positive-real-resolvent-kernel} to the end of this section, and now use this Lemma to prove Proposition \ref{prop:off-spectrum-resolvent}. We first show that
\begin{equation}\label{eq:kernel-Ls-fragment}
  \sup_{z\in\Omega_{\varepsilon}}
  \lVert K_z\rVert_{\Lp{s}}<\infty,
  \quad 1\leq s<\frac{N}{N-2}.
\end{equation}
Indeed, estimate \eqref{eq:positive-real-two-region-kernel-bound} and the weighted polar-coordinate
formula give
\begin{align*}
  \lVert K_z\rVert_{\Lp{s}}^{s}
  &\leq C_{\varepsilon,\kappa, N}
  \left(
    \int_0^1 r^{(2-N)s+N-1}\,\dd r
    +\int_1^\infty
       r^{N-1-(N-1)s/2}e^{-s\varepsilon r/2}\,\dd r
  \right).
\end{align*}
The second integral is finite for every finite $s\geq1$, owing to the
exponential decay.  The first integral is finite if and only if
\[
  (2-N)s+N>0,
\]
which is equivalent to $s<N/(N-2)$ since $N\geq3$.  Hence, this proves
\eqref{eq:kernel-Ls-fragment}.

Now choose $s>0$ according to
\[
  \frac1p+\frac1s=1+\frac1q,
\]
then the condition
$0\leq \frac1p-\frac1q<\frac2N$ ensures
$1\leq s<\frac N{N-2}$.  Since $K_z$ is radial,  we can  apply Young's inequality \eqref{ine:dunkl-young} to  obtain
\[
  \bigl\lVert(\Lk-z)^{-1}f\bigr\rVert_{\Lp{q}}
  =\lVert f*_D K_z\rVert_{\Lp{q}}
  \leq \lVert f\rVert_{\Lp{p}}
     \lVert K_z\rVert_{\Lp{s}}
  \leq C_{\varepsilon,\kappa, N, s}
     \lVert f\rVert_{\Lp{p}},
\]
uniformly for $z\in\Omega_{\varepsilon}$.

On the Sobolev line
\[
  \frac1p-\frac1q=\frac2N,
\]
the exponent required  is $s=N/(N-2)$ by Young's inequality.  This is precisely
the critical exponent at which the singular part of the kernel ceases to be
strongly integrable.  Indeed,
\[
  \int_0^1
    \bigl(r^{2-N}\bigr)^{N/(N-2)}r^{N-1}\,\dd r
  =\int_0^1\frac{\dd r}{r}
  =\infty.
\]
Consequently, $|\cdot|^{2-N}\notin\Lp{N/(N-2)}$ locally at the origin, and Young's inequality cannot be applied directly at the Sobolev
endpoint.  We instead use the Dunkl Hardy--Littlewood--Sobolev inequality.

For any radial function $h$, it follows from \eqref{radial} that
\begin{equation}\label{eq:radial-absolute-value-fragment}
  \bigl|\tau_xh(y)\bigr|
  \leq \tau_x|h|(y),\quad\forall x,y\in\mathbb{R}^d.
\end{equation}
In addition, the pointwise order is preserved for radial functions: if $0\leq h_1\leq h_2$, then
\begin{equation}\label{eq:radial-monotonicity-fragment}
  0\leq \tau_xh_1(y)\leq \tau_xh_2(y),\quad\forall x,y\in\mathbb{R}^d.
\end{equation}
Since $K_z$ is radial, using commutativity of the Dunkl convolution, \eqref{eq:positive-real-global-kernel-bound},
\eqref{eq:radial-absolute-value-fragment} and
\eqref{eq:radial-monotonicity-fragment}, we conclude that
\begin{align*}
  |f*_D K_z(x)|
  &=|K_z*_D f(x)|\\
  &=\left|
      \int_{\R^d}\tau_xK_z(-y)f(y)\omega_\kappa(y)\,\dd y
    \right|\\
  &\leq \int_{\R^d}
      \tau_x|K_z|(-y)|f(y)|\omega_\kappa(y)\,\dd y\\
  &\leq C_{\varepsilon,\kappa, N}
    \int_{\R^d}
      \tau_x\bigl(|\,\cdot\,|^{2-N}\bigr)(-y)
      |f(y)|\omega_\kappa(y)\,\dd y\\
  &\leq C_{\varepsilon,\kappa, N}I_2^{\kappa}(|f|)(x).
\end{align*}
The Dunkl Hardy--Littlewood--Sobolev theorem applied with $s=2$,
therefore yields
\[
  \|f*_D K_z\|_{\Lp{q}}
  \leq C_{\varepsilon,\kappa,N}
       \|I_2^{\kappa}(|f|)\|_{\Lp{q}}
  \leq C_{\varepsilon,\kappa,N,p}\|f\|_{\Lp{p}},
\]
whenever
\[
  1<p<\frac N2,
  \quad
  \frac1q=\frac1p-\frac2N.
\]
\subsection{Case (ii): $z\in \Omega_{\varepsilon}^{-}$.}Since $|z|=1$ and $z$ lies within distance $\varepsilon$ of the positive
half-line, we may write
$
  z=\sqrt{1-\delta^2}+i\delta,
  \  0<|\delta|\leq\varepsilon\ll1.
$
 In this case, we will prove
\begin{theorem}\label{theo:limits}
Assume that  $N\geq3$ and $(p,q)$ satisfies
\[\frac{2}{N+1}
  \leq \frac1p-\frac1q
  \leq \frac2N,\quad
  \frac{2N}{N+3}<p<\frac{2N}{N+1},
  \quad
  \frac{2N}{N-1}<q<\frac{2N}{N-3}.
\]
 Then, for $\lambda>0$, there exists a constant $C=C(\kappa,N,p,q)$ such that
\begin{equation}\label{ine:reso-re}
  \bigl\|(\Lk-(\lambda\pm i0))^{-1}f\bigr\|_{\Lp q}
  \leq C \lambda^{\frac N2(\frac1p-\frac1q)-1}\|f\|_{\Lp p},
\end{equation}
where $(\Lk-(\lambda\pm i0))^{-1}:=\lim\limits_{\epsilon\rightarrow 0^+}(\Lk-(\lambda\pm i\epsilon))^{-1}$.
\end{theorem}
We postpone the proof of the key theorem in the next section. Assuming \(\epsilon\) is sufficiently small, we can use Theorem 3.3 to obtain the estimate for the interior point \(z\), that is,
\begin{equation}\label{inner}
    \big\|(\mathcal{L}_\kappa - z)^{-1} f\big\|_{L^q_\kappa(\mathbb{R}^d)} \le C \|f\|_{L^p_\kappa(\mathbb{R}^d)},\quad\forall z\in\Omega_\varepsilon^-.
\end{equation}

Indeed, let \(\lambda = \sqrt{1 - \delta^2} > 0\) and without loss of generality, assume \(\delta > 0\). Denote the boundary limit operator as \(R_0 := (\mathcal{L}_\kappa - (\lambda + i0))^{-1}\). By Theorem 3.3, there exists a constant \(C_1 = C_1(\kappa, N, p, q) > 0\) such that
\begin{equation*}
    \|R_0\|_{L^p_\kappa \to L^q_\kappa} \le C_1 \lambda^{\frac{N}{2}(\frac{1}{p} - \frac{1}{q}) - 1}.
\end{equation*}
Since \(0 < \delta \le \epsilon \ll 1\), we have \(\lambda = \sqrt{1-\delta^2} \sim 1\), and thus \(\|R_0\|_{L^p_\kappa \to L^q_\kappa} \le \tilde{C}_1\) for some constant \(\tilde{C}_1\). By the first resolvent identity, for the interior point \(z = \lambda + i\delta\), we have
\begin{equation*}
    (\mathcal{L}_\kappa - z)^{-1} = R_0 \big(I - i\delta R_0\big)^{-1}.
\end{equation*}
Taking the operator norm and applying the Neumann series expansion, provided that \(\delta \|R_0\|_{L^p_\kappa \to L^q_\kappa} < 1\), we get
\begin{equation*}
    \|(\mathcal{L}_\kappa - z)^{-1}\|_{L^p_\kappa \to L^q_\kappa} \le \frac{\|R_0\|_{L^p_\kappa \to L^q_\kappa}}{1 - \delta \|R_0\|_{L^p_\kappa \to L^q_\kappa}}.
\end{equation*}
By taking \(\epsilon\) small enough such that \(\delta \|R_0\|_{L^p_\kappa \to L^q_\kappa} \le 1/2\), we can absorb the denominator into a new constant \(C = 2\tilde{C}_1\). Consequently, we obtain the uniform bound \eqref{inner}.

\subsection{The proof of Lemma \ref{lem:positive-real-resolvent-kernel}}

Recall that  $N=d+2\gamma_{\kappa}>3$,  $\nu=N/2-1$, and $z\in\Omega_{\varepsilon}$.  Put
\[
  m_z(\rho)=\frac{1}{\rho^2-z},
  \quad
  K_z=\Fk^{-1}m_z,
  \quad
  \rho=|\xi|.
\]
Since $m_z$ is radial, its inverse Dunkl transform $K_z$ is radial as
well.  Moreover,
\begin{equation*}
  (\Lk-z)^{-1}f
  =\Fk^{-1}\!\left(m_z\Fk f\right)
  =f*_D K_z.
\end{equation*}
We first identify \eqref{eq:positive-real-modified-bessel-kernel}.
The radial Dunkl inversion formula \eqref{eq:radialDunkl} gives
\begin{equation}
  K_z(r)
  =C_{\kappa,N}
    \int_0^\infty
      \frac{\mathfrak{j}_{\nu}(r\rho)\rho^{N-1}}
           {\rho^2-z}\dd\rho.
  \label{eq:radial-bessel-resolvent}
\end{equation}
Choose the branch
\[
  w=\sqrt{-z},
  \quad
  \operatorname{Re}w>0.
\]

\textbf{Case (i): $\operatorname{Re}z<0$.}
In this case, $\operatorname{Re}(\rho^2-z)>0$ for every $\rho\geq0$,
and therefore
\begin{equation}
  \frac{1}{\rho^2-z}
  =\int_0^\infty e^{-t(\rho^2-z)}\dd t.
  \label{equ:integral}
\end{equation}
The inverse Dunkl transform of the radial Gaussian is
\begin{equation}
  \Fk^{-1}\!\left(e^{-t|\cdot|^2}\right)(x)
  =C_{\kappa,N}t^{-N/2}e^{-|x|^2/(4t)},
  \quad t>0,
  \label{equ:radial-Gaussian}
\end{equation}
Consequently, we insert \eqref{equ:integral} and \eqref{equ:radial-Gaussian} into \eqref{eq:radial-bessel-resolvent} to obtain
\begin{align*}
  K_z(r)
  &=\int_0^\infty e^{tz}
    \Fk^{-1}\!\left(e^{-t|\cdot|^2}\right)(r)\dd t \notag\\
  &=C_{\kappa,N}\int_0^\infty
    t^{-N/2}e^{-w^2t-r^2/(4t)}\dd t.
  \label{eq:heat-kernel-resolvent}
\end{align*}
Since $N/2=\nu+1$, the Schlaefli integral representation
\cite[Eq.~(10.32.10)]{DLMF} gives
\[
  \int_0^\infty
    t^{-\nu-1}e^{-w^2t-r^2/(4t)}\dd t
  =2^{\nu+1}\left(\frac{w}{r}\right)^{\nu}\mathcal{K}_{\nu}(wr).
\]

\vspace{0.2cm}
\textbf{Case (ii): $\operatorname{Re}z>0$.}
For $\operatorname{Re}z>0$, the usual heat semigroup representation \eqref{equ:integral}
fails because
$\int_0^\infty e^{-t(\rho^2-z)}\dd t$ diverges whenever
$\rho^2<\operatorname{Re}z$.   Instead, we use the rotated
Laplace representation
\begin{equation}
  \frac{1}{\rho^{2}+w^{2}}
  =\frac{1}{w}\int_{0}^{\infty}
      e^{-s(w+\rho^{2}/w)}\,\dd s.
  \label{eq:rotated-laplace-representation}
\end{equation}
To justify the ensuing interchange of integrations, for $\eta>0$ set
\[
  m_{z,\eta}(\rho)
  :=\frac{e^{-\eta\rho^2/w}}{\rho^2+w^2},
  \quad
  K_{z,\eta}:=\Fk^{-1}m_{z,\eta}.
\]
Multiplying \eqref{eq:rotated-laplace-representation} by
$e^{-\eta\rho^2/w}$ and making the change of variables $u=s+\eta$
gives
\begin{align}
  m_{z,\eta}(\rho)
  &=\frac1w\int_0^\infty
    e^{-sw}e^{-(s+\eta)\rho^2/w}\dd s \notag\\
  &=\frac{e^{\eta w}}{w}\int_\eta^\infty
    e^{-uw}e^{-u\rho^2/w}\dd u.
  \label{eq:regularized-rotated-laplace}
\end{align}
Using the radial Dunkl inversion formula \eqref{eq:radialDunkl} and
\eqref{eq:regularized-rotated-laplace}, we obtain
\begin{align}
  K_{z,\eta}(r)
  &=C_{\kappa,N}r^{-\nu}\frac{e^{\eta w}}{w}
    \int_{\eta}^{\infty}e^{-uw}
      \left(
        \int_{0}^{\infty}e^{-u\rho^{2}/w}
          \rho^{\nu+1}J_{\nu}(r\rho)\,d\rho
      \right)du.
  \label{eq:regularized-radial-inversion}
\end{align}
For $\operatorname{Re}\beta>0$ and $\mu>-1$, the Gaussian--Bessel identity
\begin{equation*}
  \int_{0}^{\infty}e^{-\beta\rho^{2}}
    \rho^{\mu+1}J_{\mu}(r\rho)\,d\rho
  =\frac{r^{\mu}}{(2\beta)^{\mu+1}}
     e^{-r^{2}/(4\beta)}
  \label{eq:gaussian-bessel-identity}
\end{equation*}
holds, see \cite[Eq.~(10.22.51)]{DLMF}.  Substitution into
\eqref{eq:regularized-radial-inversion} yields
\begin{equation}
  K_{z,\eta}(r)
  =C_{\kappa,N}w^{\nu}e^{\eta w}
    \int_{\eta}^{\infty}u^{-\nu-1}
      e^{-w\left(u+\frac{r^{2}}{4u}\right)}du.
  \label{eq:regularized-kernel-u-integral}
\end{equation}
Let $c=\operatorname{Re}w>0$.  For every fixed $r>0$,
\[
  \left|
    u^{-\nu-1}e^{-w(u+r^2/(4u))}
  \right|
  =u^{-\nu-1}e^{-c(u+r^2/(4u))},
\]
which is integrable on $(0,\infty)$. In fact, the term
$e^{-cr^2/(4u)}$ controls the integral at the origin, while $e^{-cu}$
controls it at infinity.
Furthermore,
$m_{z,\eta}(\rho)\longrightarrow(\rho^2+w^2)^{-1}
$ as $\eta\to0$.
Hence, dominated convergence gives, for every $r>0$,
\begin{equation}
  K_z(r)
  =C_{\kappa,N}w^{\nu}
    \int_{0}^{\infty}u^{-\nu-1}
      e^{-w\left(u+\frac{r^{2}}{4u}\right)}du.
  \label{eq:unregularized-kernel-u-integral}
\end{equation}
Finally, setting $u=(r/2)e^v$ in
\eqref{eq:unregularized-kernel-u-integral} gives
\begin{align*}
  &\int_0^\infty u^{-\nu-1}
    e^{-w(u+r^2/(4u))}\dd u \\
  &\qquad
  =\left(\frac2r\right)^\nu
    \int_{-\infty}^{\infty}
      e^{-wr\cosh v}e^{-\nu v}\dd v \\
  &\qquad
  =2\left(\frac2r\right)^\nu
    \int_0^\infty
      e^{-wr\cosh v}\cosh(\nu v)\dd v \\
  &\qquad
  =2\left(\frac2r\right)^\nu \mathcal{K}_\nu(wr),
\end{align*}
where the last equality follows from
\eqref{second}.  Substituting this identity into
\eqref{eq:unregularized-kernel-u-integral} proves
\begin{equation*}
  K_z(r)
  =C_{\kappa,N}
    \left(\frac{w}{r}\right)^\nu \mathcal{K}_\nu(wr).
\end{equation*}


\vspace{0.2cm}

We next estimate \eqref{eq:positive-real-two-region-kernel-bound} and \eqref{eq:positive-real-global-kernel-bound}.  Since
$z\in\Omega_{\varepsilon}$ and $|w|=1$, then for $0<\varepsilon\ll0$, we have
\begin{equation}
  \operatorname{Re} w\geq \varepsilon/2>0.
  \label{eq:positive-real-part}
\end{equation}
Since $\nu>0$, the small-argument expansion \cite[Section~10.30]{DLMF} gives
\[
  \mathcal{K}_{\nu}(s)=2^{\nu-1}\Gamma(\nu)s^{-\nu}(1+o(1)),
  \qquad s\to0.
\]
Consequently, using \eqref{eq:positive-real-modified-bessel-kernel} and $|w|=1$, we get
\begin{equation}
  |K_z(r)|\leq C_{\varepsilon,\kappa,N}r^{-2\nu}
  =C_{\varepsilon,\kappa,N}r^{2-N},
  \qquad 0<r\leq1.
  \label{eq:near-field-bound}
\end{equation}
For $s\geq1$, the large-argument expansion \cite[Section~10.40]{DLMF} yields
\[
  \mathcal{K}_{\nu}(s)
  =\left(\frac{\pi}{2s}\right)^{1/2}e^{-s}
    \bigl(1+O(|s|^{-1})\bigr).
\]
 It follows from
\eqref{eq:positive-real-part} that
\begin{equation}
  |K_z(r)|
  \leq C_{\varepsilon,\kappa,N}
       r^{-\nu-1/2}e^{-\varepsilon r/2}
  =C_{\varepsilon,\kappa,N}
       r^{-(N-1)/2}e^{-\varepsilon r/2},\quad r\geq1.
  \label{eq:far-field-bound}
\end{equation}
Combining \eqref{eq:near-field-bound} and \eqref{eq:far-field-bound} proves \eqref{eq:positive-real-two-region-kernel-bound}.

\vspace{0.2cm}
It remains to derive the global estimate.  The claim follows directly
from \eqref{eq:near-field-bound} when $0<r\leq1$.  If $r\geq1$, then
\[
  r^{-(N-1)/2}e^{-\varepsilon r/2}
  =r^{2-N}
    \left(r^{(N-3)/2}e^{-\varepsilon r/2}\right)
  \leq C_{\varepsilon,N}r^{2-N},
\]
because $r^{(N-3)/2}e^{-\varepsilon r/2}$ is bounded on
$[1,\infty)$.  Therefore, we have proved \eqref{eq:positive-real-global-kernel-bound}.


\section{The proof of Theorem \ref{theo:limits}}In this section, we prove the critical theorem about the limitation absorbing theorem, that is, Theorem \ref{theo:limits}. By the scaling invariant of the operator $\Lk$ again, it suffices to prove \eqref{ine:reso-re} when $\lambda=1$, that is
 \begin{equation*}
  \bigl\|(\Lk-(1\pm i0))^{-1}f\bigr\|_{\Lp q}
  \leq C\|f\|_{\Lp p}.
\end{equation*}
We first have
 \begin{equation*}
  (\Lk-(1\pm i0))^{-1}f(x)
  =K^{\pm}*_D f(x)  =f*_D K^{\pm}(x).
\end{equation*}
Here $K^{\pm}$ denotes the outgoing or incoming  kernel
$$ K^{\pm}(r)
 :=\Fk^{-1}
 \left(\frac{1}{|\cdot|^{2}-(1\pm i0)}\right),\quad r=\abs{x}.$$
The radial Dunkl inversion formula yields, for \(r=\abs{x}>0\),
\begin{align*}
 K^{\pm}(r)
 =C_{\kappa,N}
   \int_{0}^{\infty}
   \frac{\mathfrak j_{\nu}(r\rho)\rho^{N-1}}
        {\rho^{2}-(1\pm i0)}\dd\rho
 =C_{\kappa,N}r^{-\nu}
   \int_{0}^{\infty}
   \frac{J_{\nu}(r\rho)\rho^{\nu+1}}
        {\rho^{2}-(1\pm i0)}\dd\rho .
\end{align*}
Then the modified Bessel representation of
the resolvent kernel  gives
\begin{equation}
 K^{\pm}(r)=c^{\pm}_{\kappa,N}r^{-\nu}H_{\nu}^{\pm}(r),
 \label{eq:boundary-hankel-kernel}
\end{equation}
where $H_{\nu}^{\pm}$  are the standard Hankel functions  \cite[Eq.~(10.27.8)]{DLMF}. Since \(\nu>0\), the small-argument
expansions of the Hankel functions \cite[Sec.~10.7(i)]{DLMF} imply
\begin{equation}
 \abs{H_{\nu}^{\pm}(r)}\leq C_{\nu}r^{-\nu},
 \quad 0<r\leq1.
 \label{pro:Hankl-short}
\end{equation}
On the other hand, the large-argument Hankel expansions
\cite[Eqs.~(10.17.5)--(10.17.6)]{DLMF} give
\begin{equation}
 H_{\nu}^{\pm}(r)
 =e^{\pm ir}r^{-1/2}a_{\pm}(r),
 \quad r\geq1,
 \label{pro:Hankl-lang}
\end{equation}
where \(a_{\pm}\in C^{\infty}[1,\infty)\) are symbols of order zero such that
\begin{equation}
 \abs{\partial_{r}^{\ell}a_{\pm}(r)}
 \leq C_{\ell}r^{-\ell},
 \quad r\geq1,\quad \ell=0,1,2,\ldots .
 \label{eq:symbol-a}
\end{equation}
Combining \eqref{eq:boundary-hankel-kernel} with
\eqref{pro:Hankl-short} gives the local bound
\begin{equation}
 \abs{K^{\pm}(r)}\lesssim r^{-2\nu}=r^{2-N},
 \quad 0<r\leq1,
 \label{eq:kernel-local}
\end{equation}
whereas \eqref{pro:Hankl-lang} gives the far-field representation
\begin{equation*}
 K^{\pm}(r)
 =e^{\pm ir}r^{-(N-1)/2}a_{\pm}(r),
 \quad r\geq1,
 \label{eq:kernel-far}
\end{equation*}
with \(a_{\pm}\) satisfying \eqref{eq:symbol-a}.

\vspace{0.2cm}
Choose a nonincreasing function \(\chi\in C_{c}^{\infty}[0,\infty)\)
such that
\[
 0\leq\chi\leq1,
 \quad
 \chi(r)=1\quad(0\leq r\leq8/5),
 \quad
 \chi(r)=0\quad(r\geq2).
\]
Define the local kernel by
\[
 K_{\mathrm{loc}}^{\pm}(r):=\chi(r)K^{\pm}(r).
\]
\begin{proposition}[local kernel bound]
\label{prop:off-spectrum-resolvent-1}Suppose that $N\geq3$ and  $(p,q)$ in Theorem \ref{theo:limits}.
Then, for every $f\in{\Lp p}$, there holds
\begin{equation}\label{eq:local-kernel-Lp-Lq}
  \bigl\|K^\pm_{loc}*_Df\bigr\|_{\Lp q}
  \lesssim \|f\|_{\Lp p}.
\end{equation}
\end{proposition}

\begin{proof} The proof is similar to Lemma  \ref{lem:positive-real-resolvent-kernel}. We give a brief proof here.
Using the polar-coordinate formula \eqref{eq:polar}, one has
\begin{align*}
\norm{K^\pm_{\mathrm{loc}}(r)}{L^s_\kappa}^s
\leq C_{\kappa,N}\int_0^2r^{(2-N)s+N-1}dr<\infty.
\end{align*}
The inequality holds because the condition $\frac1p-\frac1q<\frac2N$ ensures \(1\leq s<\frac N{N-2}\).
Young's inequality \eqref{ine:dunkl-young} therefore yields
\begin{align*}
\norm{K^\pm_{\mathrm{loc}}*_D f}{L^q_\kappa}
\le C \norm{f}{L^p_\kappa}\norm{K^\pm_{\mathrm{loc}}}{L^s_\kappa}.
\end{align*}
It remains to consider the Sobolev line
 $\frac1p-\frac1q=\frac2N$.  As  above, we utilize estimate
\eqref{eq:kernel-local} to obtain
\begin{align*}
  | K_{loc}^\pm *_D f(x)|\leq(|K_{loc}^\pm| *_D |f|(x))\lesssim (|\cdot|^{2-N}*_D|f|)(x)\lesssim I_2^\kappa(|f|)(x).
\end{align*}
Hence, Lemma \ref{lem:dunkl-HLS} gives
$$\norm{K_{\mathrm{loc}}*_D f}{L^q_\kappa}\lesssim \norm{I_2^\kappa(|f|)}{L^q_\kappa}\lesssim \norm{ f}{L^p_\kappa}.$$
\end{proof}

 \begin{proposition}\label{prop:day}
 Suppose that $N\geq3$ and  $(p,q)$ in Theorem \ref{theo:limits}.
Then, for every $f\in{\Lp p}$,
\begin{equation}\label{eq:local-kernel-Lp-Lq-2}
  \bigl\|K^\pm_{osc}*_Df\bigr\|_{\Lp q}
  \leq C\|f\|_{\Lp p},
\end{equation}
where $K^\pm_{osc}$ denotes the oscillatory part of the boundary resolvent kernel given by
\begin{equation*}K_{\mathrm{osc}}^\pm(r):=(1-\chi(r))K^{\pm}(r).\end{equation*}
 \end{proposition}
 \begin{proof}

Choose $\beta\in C_c^\infty(3/4,8/3)$ such that
$1=\sum\limits_{j\in\mathbb Z}
\beta(2^{-j}r).
$
Because \(K_{\mathrm{osc}}^{\pm}(r)=0\) for \(0<r\leq8/5\), all terms with
\(j\leq-1\) vanish on its support.  Hence
\begin{equation*}\label{eq:kjdef}
 K_{osc}^{\pm}(r)=\sum_{j\geq0}k_{j}^{\pm}(r),
 \quad
 k_{j}^{\pm}(r):=\beta(2^{-j}r)K_{osc}^{\pm}(r).
\end{equation*}
Set $Tf:=
 f*_DK_{\mathrm{osc}}^{\pm}$, then
 \begin{equation*} Tf
=\sum_{j\geq0}T_{j}^{\pm}f,\qquad T_{j}^{\pm}f:=f*_D k_{j}^{\pm},
 \quad j\geq1.
 \end{equation*}
Then we have the following estimation of the term $T_j f$.
\begin{lemma}\label{lem:dyadic}For every $j\geq0$, there holds
\begin{equation}\label{eq:Btarget}
\norm{T_jf}{L^{q_0}_\kappa}
\lesssim 2^{\frac j2}\norm{f}{L^2_\kappa},
\quad
q_0=\frac{2(N+1)}{N-1}.
\end{equation}
\end{lemma}
The proof of Lemma 4.3 is postponed to the next section. 
Recall that, for \(j\geq0\),
\begin{equation*}\label{def:kj}
 k_{j}^{\pm}(r)
 =\beta(2^{-j}r)e^{\pm i r}r^{-(N-1)/2}a_{\pm}(r).
\end{equation*}
Since \(k_{j}^{\pm}(r)\) is supported where \(r\sim 2^{j}\), there holds
\[
\norm{k_j}{L^\infty_\kappa}
\lesssim 2^{-j(N-1)/2}.
\]
Therefore, by Young's inequality \eqref{ine:dunkl-young}, it  implies
\begin{equation}\label{eq:1inf}
\norm{T_jf}{L^\infty_\kappa}
\lesssim 2^{-j(N-1)/2}\norm{f}{L^1_\kappa}.
\end{equation}
Interpolating \eqref{eq:Btarget} yields with \eqref{eq:1inf}
\begin{equation}\label{eq:dyadicpq}
\norm{T_jf}{L^q_\kappa}
\lesssim
2^{j(\frac{N+1}{2}-\frac Np)}\norm{f}{L^p_\kappa},\quad q=\frac{N+1}{N-1}p^\prime\;, 1\leq p\leq2.
\end{equation}
More generally, the support and  properties of
\(k_{j}\) imply
\[
\norm{k_j}{L^{p'}_\kappa}\lesssim
 2^{j\left(\frac{N}{p'}-\frac{N-1}{2}\right)}
=2^{j(\frac{N+1}{2}-\frac Np)}.
\]
Therefore,
\begin{equation}\label{eq:pinf}
\norm{T_jf}{L^\infty_\kappa}
\lesssim
2^{j(\frac{N+1}{2}-\frac Np)}\norm{f}{L^p_\kappa}.
\end{equation}
Interpolating further  \eqref{eq:dyadicpq} with \eqref{eq:pinf} at fixed $p$ yields
\begin{equation*}\label{eq:regionone}
\norm{T_jf}{L^q_\kappa}
\lesssim
2^{j(\frac{N+1}{2}-\frac Np)}\norm{f}{L^p_\kappa},\quad
q\ge\frac{N+1}{N-1}p'\;, 1\leq p\leq2.
\end{equation*}
The dyadic sum is absolutely convergent provided
\begin{equation*}\label{eq:psum}
\frac{N+1}{2}-\frac Np<0,
\quad\text{i.e.}\quad
p<\frac{2N}{N+1}.
\end{equation*}

\vspace{0.2cm}
 Now, we aim to prove Proposition \ref{prop:day}. We divide the proof into several cases.

{\bf Case 1:} $(\frac1p,\frac1q)$ is in the region ${\emph ABE}$ of Figure 1. We deduce that for $1\leq p<\frac{2N}{N+1}$, $ q\geq\frac{N+1}{N-1}p^\prime,$
\begin{equation}\label{eq:Btarget-2-1}
\norm{T f}{L^{q}_\kappa(\R^d)}
\lesssim \sum _{j\geq0}\norm{T_jf}{L^{q}_\kappa(\R^d)}\lesssim\norm{f}{L^p_\kappa(\R^d)}.
\end{equation}

{\bf Case 2:} $(\frac1p,\frac1q)$ is in the region ${\emph A^\prime B^\prime E^\prime}$ of Figure 1. By \eqref{eq:Btarget-2-1} and duality, we obtain
\begin{equation*}\label{eq:Btarget-2-2}
\norm{T f}{L^{q}_\kappa(\R^d)}
\lesssim \norm{f}{L^p_\kappa(\R^d)},
\end{equation*}
provided that $1\leq q^\prime<\frac{2N}{N+1}$, $ p^\prime\geq\frac{N+1}{N-1}q.$

{\bf Case 3:} $(\frac1p,\frac1q)$ is in the region ${\emph EBB^\prime E^\prime }$ of Figure 1. By Case 1 and Case 2, we use the interpolation to get \eqref{eq:local-kernel-Lp-Lq-2} for $(\frac1p,\frac1q)$ in the region ${\emph EBB^\prime E^\prime }$ of Figure 1.

{\bf Case 4:} $(\frac1p,\frac1q)$ is in the region ${\emph BB^\prime }$ of Figure 1, i.e. $\frac1p-\frac1q=\frac2{N+1}$. We follow the argument of \cite[Section~4]{FanelliZhangZheng2023}. By the real interpolation \cite[Chapter~3]{BerghLofstrom1976}, it suffices to prove
\begin{equation*}\label{eq:Btarget-2-3}
\norm{T f}{L^{q,\infty}_\kappa(\R^d)}
\leq C\norm{f}{L^{p,1}_\kappa(\R^d)},\quad (q,p)=(q_+,p_+) \ \text{or}(q_-,p_-),
\end{equation*}
where $L_\kappa^{q,\infty}$ and $L_\kappa^{p,1}$ are the Dunkl Lorentz spaces, and
\begin{align*}
 p_-&=\frac{2N(N+1)}{N^2+4N-1},
 &q_-&=\frac{2N}{N-1},                                  \label{index:pq-}\\
 p_+&=\frac{2N}{N+1},
 &q_+&=\frac{2N(N+1)}{(N-1)^2}.
\end{align*}
Since $T$ is linear, we only need to establish the corresponding restricted weak-type estimates
\begin{equation}\label{ine:Lorentz}
\sigma|\{x\in\R^d:|T\mathbf1_E(x)|>\sigma|\}|_\kappa^{\frac 1q}\leq C|E|_\kappa^{\frac1p},\quad \sigma>0,\; |E|_\kappa<\infty
\end{equation}
for measurable sets \(E\subset\R^d\) of finite measure and  $\sigma>0$, where
\[
 |E|_\kappa:=\mu_\kappa(E)
 =\int_E\omega_\kappa(x)\,dx,
\]
and $\mathbf1_E$ denotes its characteristic function. Let
\[
 A:=\{x\in\R^d:|T\mathbf1_E(x)|>\sigma\},
\]
Thus we have
\begin{equation}\label{Ine:partion}
\begin{split}
|A|_\kappa=&\int_A\omega_\kappa(x) \,dx\leq\frac 1\sigma\int_{\R^d}|T\mathbf1_E(x)|\mathbf1_A(x)\omega_\kappa(x)\,dx\\
&\leq\frac 1\sigma\sum_{j\geq0}\int_{\R^d}
       |T_j\mathbf1_E(x)|1_A(x)\omega_\kappa(x)\,dx\\
       &\leq\frac 1\sigma\Big(\sum_{0\leq j\leq J}\norm{T_j\mathbf1_E(x)}{L^{2}_\kappa(\R^d)}|A|_\kappa^{\frac12}+\sum_{ j\geq J+1}\norm{T_j\mathbf1_E(x)}{L^{\infty}_\kappa(\R^d)}|A|_\kappa\Big),
\end{split}
\end{equation}
where $J\in {N}$ will be determined lately. From  \eqref{eq:Btarget} and \eqref{eq:1inf}, we have
$$\norm{T_j}{L^1_\kappa\to L^\infty_\kappa}
\lesssim 2^{-j(N-1)/2},\quad \norm{T_j}{L^2_\kappa\to L^{q_0}_\kappa}
\lesssim 2^{j/2}.$$
By duality, we also have
$$\norm{T_j}{L^{q_0^\prime}_\kappa\to L^{2}_\kappa}
\lesssim 2^{j/2}.$$
We plug these into \eqref{Ine:partion} to show
\begin{equation*}\label{ine:case4-pq-}
\begin{split}
 \sigma|A|_\kappa
 &\leq C\Big(\sum_{0\leq j\leq J}2^{j/2}|E|_\kappa^{1/q_0'}|A|_\kappa^{\frac12}+\sum_{ j\geq J+1}2^{-j(N-1)/2}|E|_\kappa||A|_\kappa\Big)\\
 &\leq C\left(
 2^{J/2}|E|_\kappa^{1/q_0'}|A|_\kappa^{1/2}
 +2^{-J(N-1)/2}|E|_\kappa|A|_\kappa
 \right) .
 \end{split}
\end{equation*}
When $|E|_\kappa^{1/q_0}|A|_\kappa^{1/2}\geq1$, we choose $J$ such that $$2^{J/2}\leq|E|_\kappa^{1/(q_0 N)}|A|_\kappa^{1/(2 N)}\leq2^{(J+1)/2}.$$ Consequently
\begin{equation*}
\sigma|A|_\kappa\leq C|E|_\kappa^{\frac1{q_0^\prime}+\frac1{q_0 N}}|A|_\kappa^{\frac12+\frac1{2N}}=C|E|_\kappa^{\frac1{p_-}}|A|_\kappa^{1-\frac1{q_-}}.
\end{equation*}
This implies
\begin{equation}\label{aim1}
 \sigma|A|_\kappa^{1/q_-}
 \leq C|E|_\kappa^{1/p_-}.
\end{equation}
When $|E|_\kappa^{1/q_0}|A|_\kappa^{1/2}<1$, using \eqref{eq:1inf},  we get
\begin{equation*}
\begin{split}
 \sigma|A|_\kappa
&\lesssim \sum_{j\geq0} 2^{-\frac{N-1}{2}j }|E|_\kappa|A|_\kappa\\
&\lesssim|E|_\kappa|A|_\kappa\\
&\leq |E|_\kappa|A|_\kappa (|E|_\kappa^{\frac1{q_0}}|A|_\kappa^{\frac12})^{-\frac{N-1}N} \\
&\leq C|E|_\kappa^{\frac1{q_0^\prime}+\frac1{q_0 N}}|A|_\kappa^{\frac12+\frac1{2N}},
\end{split}
\end{equation*}
which also gives \eqref{aim1}.  To prove \eqref{ine:Lorentz} with $(p,q)=(p_+,q_+)$, similarly we have
\begin{equation}\label{ine:partion}
\begin{split}
|A|_\kappa
       &\lesssim\frac 1\sigma\Big(\sum_{0\leq j\leq J}\norm{T_j\mathbf1_E(x)}{L^{q_0}_\kappa(\R^d)}|A|_\kappa^{\frac1{q_0^\prime}}+\sum_{ j\geq J+1}\norm{T_j\mathbf1_E(x)}{L^{\infty}_\kappa(\R^d)}|A|_\kappa\Big)\\
       &\lesssim\frac 1\sigma\Big(\sum_{0\leq j\leq J}2^{j/2}|E|_\kappa^{1/2}|A|_\kappa^{\frac1{q_0^\prime}}+\sum_{ j\geq J+1}2^{-j(N-1)/2}|E|_\kappa||A|_\kappa\Big)\\
 &\lesssim\frac 1\sigma\left(
 2^{J/2}|E|_\kappa^{1/2}|A|_\kappa^{1/q_0^\prime}
 +2^{-J(N-1)/2}|E|_\kappa|A|_\kappa
 \right).
\end{split}
\end{equation}
When $|A|_\kappa^{1/q_0}|E|_\kappa^{1/2}\geq1$, we choose $J$ such that $$2^{J/2}\leq|A|_\kappa^{1/(q_0 N)}|E|_\kappa^{1/(2 N)}\leq2^{(J+1)/2}.$$ Therefore
\begin{equation*}
\sigma|A|_\kappa \leq C|A|_\kappa^{\frac1{q_0^\prime}+\frac1{q_0 N}}|E|_\kappa^{\frac12+\frac1{2N}}=C|E|_\kappa^{\frac1{p_+}}|A|_\kappa^{1-\frac1{q_+}}.
\end{equation*}
This yields
\begin{equation}\label{aim2}
 \sigma|A|_\kappa^{1/q_+}
 \leq C|E|_\kappa^{1/p_+}.
\end{equation}
When $|A|_\kappa^{1/q_0}|E|_\kappa^{1/2}<1$, using \eqref{eq:1inf},  we get
\begin{equation*}
\begin{split}
 \sigma|A|_\kappa
 \leq \sum_{j\geq0} 2^{-\frac{N-1}{2}j}|E|_\kappa|A|_\kappa\leq C|A|_\kappa^{\frac1{q_0^\prime}+\frac1{q_0 N}}|E|_\kappa^{\frac12+\frac1{2N}},
\end{split}
\end{equation*}
which means \eqref{aim2}.
In sum, we conclude the proof of the estimate \eqref{eq:local-kernel-Lp-Lq-2}.

\end{proof}

\section{The proof of Lemma \ref{lem:dyadic}} The purpose of this section is to finish the proof of  Lemma \ref{lem:dyadic}. We restrict attention to the plus sign  $e^{ir}$, the case $e^{-ir}$ can be handled in an entirely analogous manner.
We first establish the following estimate for the radial multiplier associated with a dyadic kernel.

\begin{lemma}[Dyadic multiplier estimate]\label{lem:B1}
Let  \(R=2^j\). Define
\[
 k_j(r):=\beta(r/R)e^{ir}r^{-(N-1)/2}a_+(r),\quad j\geq1,
\]
where \(a_+\in C^\infty[1,\infty)\)
satisfies
\begin{equation*}
 \abs{\partial_r^\ell a_+(r)}\leq C_\ell r^{-\ell},
 \quad r\geq1,\; \ell\geq0.
 \label{eq:symbol-a-1}
\end{equation*}
Let
$
 m_j(\rho)=\Fk(k_j)(\rho)
$
be its radial Dunkl transform.  Then, for every integer \(M\geq1\),
\begin{equation}
 \abs{m_j(\rho)}
 \leq C_M R\bigl(1+R\abs{\rho-1}\bigr)^{-M},
 \quad \rho\geq0,
 \label{eq:B1point}
\end{equation}
where \(C_M\) is independent of \(j\) and \(\rho\).  Consequently,
by choosing \(M>N/2\),
\begin{equation}
 \left(
  \int_0^\infty\abs{m_j(\rho)}^2\rho^{N-1}\dd\rho
 \right)^{1/2}
 \leq C R^{1/2}=C2^{j/2}.
 \label{eq:B1L2}
\end{equation}
\end{lemma}

\begin{proof}

The radial Dunkl-transform formula
 gives
\begin{equation}
 m_j(\rho)
 =C_{\kappa,N}\int_0^\infty
 e^{ir}\beta(r/R)a_+(r)
 r^{(N-1)/2}\mathfrak j_\nu(r\rho)\dd r.
 \label{eq:B1start}
\end{equation}
Since the support of $\beta$ is contained in $(3/4,8/3)$, the integral in
\eqref{eq:B1start} supports in $( 3R/4,8R/3)$. In the sequel, we split the range of
$\rho$  into three cases and treat each case separately.

\medskip
\textbf{Case (i): $0\le\rho\le1/2$.}
We first consider the very low subregion $0\le\rho\le2/R$.  Define
\[
A_{j,\rho}(r):=
\beta(r/R)a_+(r)r^{(N-1)/2}\mathfrak j_\nu(r\rho).
\]
On the support of $A_{j,\rho}$, one has $0\leq r\rho\leq 16/3$. Therefore,  \eqref{eq:besselsmall} and $\rho\lesssim R^{-1}$ imply that, for every integer $\ell\geq0$
\begin{equation*}\label{eq:Aders}
\abs{\partial_r^\ell A_{j,\rho}(r)}
\le C_\ell R^{(N-1)/2-\ell}.
\end{equation*}
Integrating \eqref{eq:B1start} by parts $L$ times against $e^{ir}$, with the boundary terms vanishing by virtue of the compact support of
$\beta$,  yields
\begin{align*}
 \abs{m_j(\rho)}
 &\leq C_L
 \int_{\frac{3}4 R}^{\frac{8}3R}
 \abs{\partial_r^L A_{j,\rho}(r)}\dd r\leq C_L R^{(N+1)/2-L}.
\end{align*}
Given \(M\geq1\), choose an integer
\(L\geq M+(N+1)/2\).  Then
\begin{equation*}
 \abs{m_j(\rho)}\leq C_MR^{-M},
 \quad 0\leq\rho\leq2/R.
 \label{eq:verylow}
\end{equation*}

Next, suppose that \(2/R\leq\rho\leq1/2\).  Since
\(r\geq3R/4\) on the support of \(\beta(r/R)\), we have
\(r\rho\geq3/2\), so the decomposition
\eqref{eq:bessellarge} applies.  Moreover,
\[
 r^{(N-1)/2}(r\rho)^{-(N-1)/2}
 =\rho^{-(N-1)/2}.
\]
It follows that
\begin{equation}
 m_j(\rho)
 =C_{\kappa,N}\rho^{-(N-1)/2}
 \sum_{\pm}I_\pm(\rho),
 \label{eq:Ipmrep}
\end{equation}
where
\begin{equation*}\label{eq:Ipmdef}
I_\pm(\rho)
=\int e^{ir(1\pm\rho)}B_{\pm,j,\rho}(r)\dd r,
\quad
B_{j,\rho}^\pm(r)=\beta(r/R)a_+(r)b_\pm(r\rho).
\end{equation*}
The chain rule and \eqref{eq:bsymbol} yield
\begin{equation}
 \abs{\partial_r^\ell B_{j,\rho}^\pm(r)}
 \leq C_\ell R^{-\ell}.
 \label{eq:Bders}
\end{equation}
Indeed, each derivative falling on \(b_\pm(r\rho)\) contributes a
factor satisfying
\[
 \rho^k\abs{b_\pm^{(k)}(r\rho)}
 \lesssim \rho^k(r\rho)^{-k}
 =r^{-k}\lesssim R^{-k}.
\]
Since \(\abs{1-\rho}\geq1/2\) and \(1+\rho\geq1\), repeated
integration by parts $L$ times gives
\[
 \abs{I_\pm(\rho)}
 \leq C_L R^{1-L}.
\]
Because \(\rho\geq2/R\),
\[
 \rho^{-(N-1)/2}\lesssim R^{(N-1)/2},
\]
and consequently
\[
 \abs{m_j(\rho)}
  \leq C_L R^{(N+1)/2-L}.
\]
Choosing \(L\geq M+(N+1)/2\) once again yields
\begin{equation}
 \abs{m_j(\rho)}\leq C_MR^{-M},
 \quad 0\leq\rho\leq1/2.
 \label{eq:lowall}
\end{equation}

\medskip
\textbf{Case (ii):  $1/2\le\rho\le2$.} In this case,
the representation \eqref{eq:Ipmrep} remains valid, and $\rho^{-(N-1)/2}\sim1$.  We first estimates $I_+$. The phase derivative with respect to $r$ is $1+\rho\ge3/2$. Repeating integration by parts, we obtian
\begin{equation}\label{eq:Iplus}
\abs{I_+(\rho)}\leq C_L R^{1-L},\quad L=0,1,2,\cdots
\end{equation}
The only resonant term is
\[
I_-(\rho)=\int_0^\infty e^{ir(1-\rho)}B_{j,\rho}^-(r)\dd r.
\]
If $R\abs{1-\rho}\le1$,then the trivial estimate, the support
condition, and \eqref{eq:Bders} with \(\ell=0\) give
\[
 \abs{I_-(\rho)}\lesssim R.
\]
Since \(1+R\abs{1-\rho}\leq2\) in this case, it follows that
\begin{equation*}
 \abs{I_-(\rho)}
 \leq C_MR\bigl(1+R\abs{1-\rho}\bigr)^{-M}.
 \label{eq:Iminus-near}
\end{equation*}
If $R\abs{1-\rho}\ge1$,  integrating by parts \(M\) times gives
\[
I_-(\rho)
=(-1)^M[i(1-\rho)]^{-M}
\int_0^\infty e^{ir(1-\rho)}\partial_r^M B_{j,\rho}^-(r)\dd r.
\]
Consequently, \eqref{eq:Bders} implies
\begin{align*}
 \abs{I_-(\rho)}
 \leq C_M
 \abs{1-\rho}^{-M}R^{1-M}
 =R\bigl(R\abs{1-\rho}\bigr)^{-M}
 \leq C_MR\bigl(1+R\abs{1-\rho}\bigr)^{-M},
\end{align*}
where the final inequality holds because of \(R\abs{1-\rho}\geq1\).  Combining
the two subcases, we obtain
\begin{equation*}
 \abs{I_-(\rho)}
 \leq C_MR\bigl(1+R\abs{1-\rho}\bigr)^{-M}.
 \label{eq:Iminus}
\end{equation*}
Together with \eqref{eq:Iplus}, this proves
\begin{equation}\label{eq:mid}
\abs{m_j(\rho)}
\le C_M R(1+R\abs{\rho-1})^{-M},
\quad 1/2\le\rho\le2.
\end{equation}

\medskip
\textbf{Case (iii):  $\rho\ge2$.}
In this case,
\[
 \abs{1-\rho}=\rho-1\geq\frac{\rho}{2},
 \quad
 1+\rho\geq\rho.
\]
Hence, neither phase is stationary.  Integrating \(I_\pm\)
by parts \(M\) times and using \eqref{eq:Bders} yield
\[
 \abs{I_\pm(\rho)}
 \leq C_M \rho^{-M}R^{1-M}.
\]
Substitution into \eqref{eq:Ipmrep} gives the following estimate
\begin{equation}
 \abs{m_j(\rho)}
 \leq C_M
 R^{1-M}\rho^{-M-(N-1)/2}.
 \label{eq:highstrong}
\end{equation}
Since \(\rho\geq2\) implies
\(\abs{\rho-1}\sim\rho\), and  \(R\geq2\), we have
\[
 1+R\abs{\rho-1}\sim R\rho.
\]
Moreover, \(\rho^{-(N-1)/2}\leq1\).  Hence
\eqref{eq:highstrong} implies
\begin{equation}
 \abs{m_j(\rho)}
 \leq C_MR\bigl(1+R\abs{\rho-1}\bigr)^{-M},
 \quad \rho\geq2.
 \label{eq:high}
\end{equation}
Combining \eqref{eq:lowall},\eqref{eq:mid},\eqref{eq:high} together, we  finally obtain
\begin{equation*}
 \abs{m_j(\rho)}
 \leq C_M R\bigl(1+R\abs{\rho-1}\bigr)^{-M},
 \quad \rho\geq0.
\end{equation*}

\medskip
It remains to prove the weighted \(L^2\) estimate.  Taking
\(M>N/2\), by \eqref{eq:B1point}, we have
\begin{equation}
 \int_0^\infty\abs{m_j(\rho)}^2\rho^{N-1}\dd\rho
 \lesssim
 R^2\int_0^\infty
 \bigl(1+R\abs{\rho-1}\bigr)^{-2M}
 \rho^{N-1}\dd\rho.
 \label{eq:L2-start}
\end{equation}
On \(0\leq\rho\leq2\), the weight \(\rho^{N-1}\) is bounded.
The change of variables \(u=R(\rho-1)\) therefore gives
\begin{align*}
 R^2\int_0^2
 \bigl(1+R\abs{\rho-1}\bigr)^{-2M}
 \rho^{N-1}\dd\rho
 &\lesssim
 R\int_{-R}^{R}(1+\abs{u})^{-2M}\dd u\lesssim R.
\end{align*}
On the other hand, when \(\rho\geq2\),
\[
 \bigl(1+R\abs{\rho-1}\bigr)^{-2M}
 \lesssim_M R^{-2M}\rho^{-2M}.
\]
Since \(2M>N\), the tail is integrable and satisfies
\begin{align*}
 R^2\int_2^\infty
 \bigl(1+R\abs{\rho-1}\bigr)^{-2M}
 \rho^{N-1}\dd\rho
 &\lesssim_M
 R^{2-2M}\int_2^\infty\rho^{N-1-2M}\dd\rho\\
 &\lesssim_M R^{2-2M}
 \lesssim R,
\end{align*}
where we have used \(R\geq2\) and \(M\geq1\).
Combining the two regions in \eqref{eq:L2-start} gives
\[
 \int_0^\infty\abs{m_j(\rho)}^2\rho^{N-1}\dd\rho
 \lesssim R.
\]
which is \eqref{eq:B1L2}.

\end{proof}


 For
$f\in\mathcal S(\R^d)$, define the restriction operator
\[
  \mathcal R_{\kappa}f(\theta)
  :=\Fk f(\theta),
  \quad \theta\in\Sph^{d-1},
\]
and, for $\rho>0$, define the adjoint extension operator by
\[
  \mathcal E_{\rho}g(x)
  :=\frac{1}{c_{\kappa}}
  \int_{\Sph^{d-1}}
    g(\theta)E_{\kappa}(i\rho\theta,x)\,
    \omega_\kappa (\theta)d\sigma(\theta).
\]
\begin{lemma}
\label{lem:dunkl-stein-tomas}
Let
$
q_{0}
=
\frac{2(N+1)}{N-1}.
$
Then, for every $g\in L^2(\Sph^{d-1}, \omega_\kappa (\theta)d\sigma(\theta))$,
\begin{equation}
\label{eq:dunkl-extension-ST}
\|\mathcal E_{1}g\|_
{L_{\kappa}^{q_{0}}(\mathbb R^d)}
\le
C
\|g\|_
{L^2(\mathbb S^{d-1},\omega_\kappa (\theta)d\sigma(\theta))}.
\end{equation}
Moreover, for every $\rho>0$,
\begin{equation}\label{eq:dunkl-extension-ST-scaling}
\|\mathcal E_{\rho}g\|_
{L_{\kappa}^{q_{0}}(\mathbb R^d)}
\le
C\rho^{-N/q_0}
\|g\|_
{L^2(\mathbb S^{d-1},\omega_\kappa (\theta)d\sigma(\theta))},
\end{equation}
where $C$ is independent of $\rho$ and $g$.
\end{lemma}

\begin{proof}
We invoke the Fourier--Dunkl restriction theorem for the sphere.
Theorem~A(iii) in Senapati--Boggarapu--Mondal--Mejjaoli
\cite{SenapatiEtAl2025Restriction} states that, for the sphere
\[
S_2
=
\{(\xi,\zeta)\in\mathbb R^d\times\mathbb R^n:
|\xi|^2+|\zeta|^2=1\},
\]
one has
\[
\left(
\int_{S_2}
|\widehat f(\xi,\zeta)|^2
h_{\kappa}^2(\xi)\,d\sigma(\xi,\zeta)
\right)^{1/2}
\le
C
\|f\|_{L_{\kappa}^p(\mathbb R^d\times\mathbb R^n)}
\]
whenever
\[
1\le p\le
\frac{2(N_{\kappa}+1)}{N_{\kappa}+3},
\qquad
N_{\kappa}:=d+n+2\gamma_{\kappa}.
\]
Taking \(n=0\), we have
\[
N_{\kappa}=d+2\gamma_{\kappa}=N,
\]
and \(S_2\) reduces to the unit sphere
\(\mathbb S^{d-1}\subset\mathbb R^d\). Hence, at the endpoint
\[
q_{0}^\prime=\frac{2(N+1)}{N+3},
\]
we obtain
\[
\|\mathcal R_{\kappa}f\|_
{L^2(\mathbb S^{d-1},\omega_\kappa (\theta)d\sigma(\theta)}
\le
C
\|f\|_
{L_{\kappa}^{q_{0}^\prime}(\mathbb R^d)}.
\]
By duality argument, it is obviously to get
\[
\|\mathcal E_{1}g\|_
{L_{\kappa}^{q_{0}}(\mathbb R^d)}
\le
C
\|g\|_
{L^2(\mathbb S^{d-1},\omega_\kappa (\theta)d\sigma(\theta)},
\]
which proves \eqref{eq:dunkl-extension-ST}. Finally, the homogeneity of the Dunkl kernel gives
\[
  E_{\kappa}(i\rho\theta,x)=E_{\kappa}(i\theta,\rho x).
\]
Hence
$$\mathcal E_{\rho}g(x)=\mathcal E_1g(\rho x).$$
Because $\omega_\kappa(x)d\sigma(x)$ is homogeneous of degree $N$, the change of
variables $y=\rho x$ yields
\[
  \norm{\mathcal E_{\rho}g}{L_{\kappa}^{q_0}(\R^d)}
  =\rho^{-N/q_0}
   \norm{\mathcal E_1g}{L_{\kappa}^{q_0}(\R^d)}.
\]
This Combining  with \eqref{eq:dunkl-extension-ST} proves
\eqref{eq:dunkl-extension-ST-scaling}.

\end{proof}

The desired dyadic $L^2$
  bound follows immediately from Lemma \ref{lem:B1} and Lemma \ref{lem:dunkl-stein-tomas}.

\begin{proof}
The Dunkl convolution multiplier identity gives
\[
\Fk(T_jf)(\xi)=m_j(\abs{\xi})\Fk f(\xi).
\]
Using Dunkl inversion and the weighted  polar decomposition \eqref{eq:polar}, we obtain
\[
T_jf(x)
=C\int_0^\infty
m_j(\rho)\rho^{N-1}
\mathcal E_\rho\bigl(\Fk f(\rho\,\cdot)\bigr)(x)\dd\rho.
\]
By making use of Minkowski inequality and \eqref{eq:dunkl-extension-ST-scaling}, one has
\begin{align*}
\norm{T_jf}{L^{q_0}_\kappa}
&\le C\int_0^\infty
\abs{m_j(\rho)}
\rho^{N-1-N/q_0}
\norm{\Fk f(\rho\,\cdot)}{L^2(\mathbb S^{d-1},\omega_\kappa (\theta)d\sigma(\theta))}\dd\rho.
\end{align*}
Applying the Cauchy--Schwarz inequality with respect to
the measure $\rho^{N-1}\,d\rho$, we find
\begin{align}\label{eq:CS}
 \|T_jf\|_{L_\kappa^{q_0}}
 &\leq C
 \left(
   \int_0^\infty
   |m_j(\rho)|^2
   \rho^{N-1-2N/q_0}\,d\rho
 \right)^{1/2}
 \nonumber\\
 &\qquad\times
 \left(
   \int_0^\infty
   \|\mathcal F_\kappa f(\rho\,\cdot)\|_
   {L^2(\mathbb S^{d-1},\omega_\kappa (\theta)d\sigma(\theta))}^2
   \rho^{N-1}\,d\rho
 \right)^{1/2}.
\end{align}
By the Dunkl Plancherel theorem and the polar
decomposition, the second factor on the
right-hand side of  \eqref{eq:CS} is bounded by
$
 \|f\|_{L_\kappa^2}.$
It therefore remains to prove
\begin{equation}\label{ine:mj-est}
 \int_0^\infty
 |m_j(\rho)|^2
 \rho^{N-1-2N/q_0}\,d\rho
 \leq CR.
 \end{equation}
Set
\[
 \alpha
 :=N-1-\frac{2N}{q_0}
 =\frac{N-1}{N+1}.
\]
In particular, $0<\alpha<1$. We split the integral in
\eqref{ine:mj-est} into three regions.
\[
 I_1:=\int_0^{1/2}|m_j(\rho)|^2\rho^\alpha\,d\rho,
 \qquad
 I_2:=\int_{1/2}^{2}|m_j(\rho)|^2\rho^\alpha\,d\rho,
 \qquad
 I_3:=\int_2^\infty|m_j(\rho)|^2\rho^\alpha\,d\rho.
\]\\
On  $1/2\leq\rho\leq2$, the weights
$\rho^\alpha$ and $\rho^{N-1}$ are comparable. Hence \eqref{eq:B1L2} gives
\begin{equation*}
 I_2
 \leq C\int_{1/2}^{2}
       |m_j(\rho)|^2\rho^{N-1}\,d\rho
 \leq CR.
\label{ine:mj-est-2}
\end{equation*}
If  $0<\rho\leq1/2$, then $|\rho-1|\geq1/2$.
Consequently,  \eqref{eq:B1point} implies,
for every integer $L\geq1$,
\[
 |m_j(\rho)|
 \leq C_LR(1+R|\rho-1|)^{-L}
 \leq C_LR^{1-L}.
\]
Since $\alpha>-1$, it follows that
\begin{equation}
 I_1
 \leq C_LR^{2-2L}
       \int_0^{1/2}\rho^\alpha\,d\rho
 \leq C_LR^{2-2L}.
 \label{ine:mj-est-1}
\end{equation}\\
If $\rho\geq2$, then $|\rho-1|\geq\rho/2$, and hence
\[
 |m_j(\rho)|
 \leq C_LR(1+R|\rho-1|)^{-L}
 \leq C_LR^{1-L}\rho^{-L}.
\]
Choosing $L$ such that $2L>\alpha+1.$ Then
\begin{equation}
 I_3
 \leq C_LR^{2-2L}
       \int_2^\infty\rho^{\alpha-2L}\,d\rho
 \leq C_LR^{2-2L}.
 \label{ine:mj-est-3}
\end{equation}
Since $R=2^j\geq1$ and $L\geq1$, equations \eqref{ine:mj-est-1}- \eqref{ine:mj-est-3} imply
\[
 I_1+I_2+I_3\leq CR.
\]
This proves \eqref{ine:mj-est}. Hence, \ref{lem:dyadic} holds.
\end{proof}

\end{document}